\documentclass[12pt,reqno]{amsart}

\usepackage{amsmath, amsfonts, amssymb, amsthm, amscd, amsbsy}
\usepackage{fancyhdr}
\usepackage[usenames,dvipsnames,svgnames,x11names,hyperref]{xcolor}
\usepackage{geometry}
\usepackage{graphicx}
\usepackage[pagebackref]{hyperref}%
\usepackage{amsmath}%
\usepackage{amsfonts}%
\usepackage{amssymb}
\usepackage{enumerate}

\providecommand{\U}[1]{\protect\rule{.1in}{.1in}}
\hypersetup{backref=true,
pagebackref=true,
hyperindex=true,
colorlinks=true,
breaklinks=true,
urlcolor=NavyBlue,
linkcolor=Fuchsia,
bookmarks=true,
bookmarksopen=false,
filecolor=black,
citecolor=ForestGreen,
linkbordercolor=red
}
\allowdisplaybreaks

\newtheorem{theorem}{Theorem}[section]

\newtheorem{prop}[theorem]{Proposition}

\newtheorem{lemma}[theorem]{Lemma}
\newtheorem{remark}[theorem]{Remark}
\newtheorem{example}[theorem]{Example}
\newtheorem{examples}[theorem]{Examples}
\newtheorem{foo}[theorem]{Remarks}

\def\vint{\mathop{\mathchoice%
          {\setbox0\hbox{$\displaystyle\intop$}\kern 0.22\wd0%
           \vcenter{\hrule width 0.6\wd0}\kern -0.82\wd0}%
          {\setbox0\hbox{$\textstyle\intop$}\kern 0.2\wd0%
           \vcenter{\hrule width 0.6\wd0}\kern -0.8\wd0}%
          {\setbox0\hbox{$\scriptstyle\intop$}\kern 0.2\wd0%
           \vcenter{\hrule width 0.6\wd0}\kern -0.8\wd0}%
          {\setbox0\hbox{$\scriptscriptstyle\intop$}\kern 0.2\wd0%
           \vcenter{\hrule width 0.6\wd0}\kern -0.8\wd0}}%
          \mathopen{}\int}

\begin{document}

\title[Higher order Brezis-Nirenberg problem on hyperbolic spaces]{Higher order Brezis-Nirenberg problem on hyperbolic spaces: Existence, nonexistence and symmetry of solutions }
\date{}

\author{Jungang Li}
\address{Jungang Li: Department of Mathematics\\
Brown University\\
Providence, RI 02912, USA}
\email{jungang\_li@brown.edu}
\author{Guozhen Lu}
\address{Guozhen Lu: Department of Mathematics\\
University of Connecticut\\
Storrs, CT 06269, USA}
\email{guozhen.lu@uconn.edu}
\author{Qiaohua Yang}
\address{Qiaohua Yang: School of Mathematics and Statistics, Wuhan University, Wuhan, 430072, People's  Republic of China}
\email{qhyang.math@whu.edu.cn}

\thanks{The second author was supported by a Simons Collaboration grant from the Simons Foundation. The third  author was  supported by the National Natural Science Foundation of China(No.12071353).}

\begin{abstract}
The main purpose of this paper is to establish the existence, nonexistence and symmetry of nontrivial solutions to the higher order Brezis-Nirenberg problems associated with the GJMS operators $P_k$ on bounded domains in the hyperbolic space $\mathbb{H}^n$ and as well as on the entire hyperbolic space $\mathbb{H}^n$. Among other techniques, one of our main novelties is to use crucially the Helgason-Fourier analysis on hyperbolic spaces and the higher order Hardy-Sobolev-Maz'ya inequalities and careful study of delicate properties of Green's functions of $P_k-\lambda$ on hyperbolic spaces which are of independent interests in dealing with such problems. Such Green's functions allow us to obtain the integral representations of solutions and thus to avoid using the maximum principle to establish the symmetry of solutions.
\end{abstract}

\maketitle {\small {\bf Keywords:}  Higher order Brezis-Nirenberg problem, hyperbolic spaces, Existence, nonexistence and symmetry, Helgason-Fourier analysis, Hardy-Sobolev-Maz'ya inequalities, moving plane method in integral form on hyperbolic spaces. \\

 }

\section{Introduction}

In the celebrated work of Brezis-Nirenberg \cite{BrezisNirenberg1}, the following second order semilinear equation  has been studied:

\begin{equation}\label{Brezis-Nirenberg}
  \begin{cases}
    - \triangle u = \lambda u + u^{q - 1} \ \  &\text{on } \Omega \\
    u > 0 \ \ &\text{on } \Omega \\
    u = 0 \ \ &\text{on } \partial \Omega,
  \end{cases}
\end{equation}
where $\Omega$ is a bounded domain in $\Bbb{R}^n$ for $n\ge 3$ and $q = \frac{2n}{n - 2}$ is the critical Sobolev exponent. If we denote

$$
  \Lambda_1 ( - \triangle , \Omega) = \inf_{u \in W^{1,2}_0 (\Omega)} \frac{\int_\Omega |\nabla u|^2 dx}{\int_\Omega |u|^2 dx}
$$
as the first eigenvalue of $- \triangle$ for the Dirichlet problem on $\Omega$, it was proved in \cite{BrezisNirenberg1} that when $n \geq 4$ and $\lambda \in (0 , \Lambda_1 ( -\triangle , \Omega ))$, there exists a $W^{1,2}_0 (\Omega)$ weak solution to \eqref{Brezis-Nirenberg}. On the other hand, when $\lambda \not\in (0, \Lambda_1 ( -\triangle , \Omega ))$ and if further assume $\Omega$ to be a star-shaped domain, then there is no solution to the equation \eqref{Brezis-Nirenberg}. Moreover, a solution gap phenomenon appears in the lower dimension case $n = 3$. To be precise, in the same paper \cite{BrezisNirenberg1}, it was proved by Brezis and Nirenberg that when $\Omega$ is a ball, \eqref{Brezis-Nirenberg} is solvable in dimension 3 if and only if $\lambda \in \left( \frac{1}{4} \Lambda_1 ( -\triangle , \Omega ) , \Lambda_1 ( -\triangle , \Omega ) \right)$.
This problem has since been called the  well-known Brezis-Nirenberg problem. There have been  tremendous amount of works in related problems of Brezis-Nirenberg type over the past decades.


It is well understood that the existence and nonexistence of solutions for semilinear equations are crucially depending on their nonlinearity terms. When associated with the critical Sobolev exponent, the standard variational method fails due to the lack of compactness. The lack of compactness appears in many variational problems in PDEs and geometry. One of the most well known examples is the Yamabe problem, which asks for the existence of a solution for the following equation on the Riemannian manifold $(M,g)$:

$$
 4 \frac{n-1}{n-2} \triangle_M u + R u = \tilde{R} u^{\frac{n+2}{n-2}},
$$
where $\triangle_M$ is the Laplace-Beltrami operator, $R$ is the scalar curvature and $\tilde{R}$ is a given constant. The Yamabe problem was initiated by Yamabe \cite{Yamabe1} and completely solved after a long time of effort. For the proof, see Trudinger \cite{Trudinger1}, Aubin \cite{Aubin2} and Schoen \cite{Schoen1} and for the full list of reference of the history and application of Yamabe problem, see \cite{LeeParker1}.  for hyperbolic spaces $\Bbb{H}^n$, see \cite{Benguria1,ManciniSandeep1,Stapelkamp1}.

\medskip

One of the higher order versions of Problem \eqref{Brezis-Nirenberg} can be formulated as follows:

\begin{equation}\label{Pucci-Serrin}
  \begin{cases}
    (- \triangle)^k u = \lambda u + |u|^{q - 2} u \ \ \ &\text{on } \Omega \\
    u = D u = \cdots = D^\alpha u = 0 \ \ \ &\text{on } \partial \Omega \\
    |\alpha| \leq k-1,
  \end{cases}
\end{equation}
where $\Omega\subset \mathbb{R}^n$ is a bounded domain, $n > 2k$ and $q = \frac{2n}{n - 2k}$ is the corresponding critical Sobolev exponent. Gazzola \cite{Gazzola1} proved the following existence result:

\begin{theorem}\label{Gazzola theorem}
Denote $\Lambda_1 ( (- \triangle)^k  , \Omega)$ as the first eigenvalue of $(-\triangle)^k$ on $\Omega$, then
  \begin{itemize}
    \item If $n \geq 4k$, for for every $\lambda \in (0 , \Lambda_1( (-\triangle)^k , \Omega))$ there exists a solution $u \in W^{k,2}_0(\Omega)$ to the Dirichlet problem (\ref{Pucci-Serrin}).

    \item If $2k + 1 \leq n \leq 4k - 1$, then there exists  $0 < \overline{\Lambda} < \Lambda_1( (-\triangle)^k , \Omega)) $ such that for every $\lambda \in (\overline{\Lambda} , \Lambda_1( (-\triangle)^k , \Omega)))$ the Dirichlet problem (\ref{Pucci-Serrin}) has a solution $u \in W^{k,2}_0(\Omega)$.
  \end{itemize}
\end{theorem}

\begin{remark}
  When $\Omega = \Bbb{B}^n$, an Euclidean ball, Grunau \cite{Grunau2} established stronger results. When $n \geq 4k$, it has been shown that the solution in Theorem \ref{Gazzola theorem} belongs to $C^\infty(\Bbb{B}^n) \cap C^{2k + 1}(\overline{\Bbb{B}^n})$ and the solution is positive, radially symmetric and decreasing. When $2k + 1 \leq n \leq 4k - 1$, besides the result in Theorem \ref{Gazzola theorem}, there exists $ \underline{\Lambda}$ satisfying $0 < \underline{\Lambda} \leq \overline{\Lambda} < \Lambda_1( (-\triangle)^k , \Omega))$ such that for every $\lambda \in (0 , \underline{\Lambda})$, the Dirichlet problem (\ref{Pucci-Serrin}) has no positive solution.
\end{remark}

 On the other hand, if $\Omega$ is a star-shaped domain and $k \geq 2$, it was shown in \cite{PucciSerrin2} that the Problem \eqref{Pucci-Serrin} has no solution when $\lambda < 0$ (see also Chapter 7 of \cite{GazzolaGrunauSweers1} for a full list of references for related results).

\medskip

In lower dimension cases $n = 2k + 1, 2k + 2, \cdots, 4k - 1$, a similar solution gap phenomenon has been observed as in the second order equation case, i.e. there exits a positive constant $\underline{\Lambda} ( (-\triangle)^k , \Omega )$ such that there exists a solution $u \in W^{k,2}_0 (\Omega)$  to the Problem \eqref{Pucci-Serrin} for every $\lambda \in ( \underline{\Lambda} ((-\triangle)^k , \Omega) , \Lambda_1 ( (- \triangle)^k  , \Omega))$. When $\Omega$ is a ball, Pucci and Serrin \cite{PucciSerrin1} further conjectured that for dimensions $n = 2k + 1, 2k + 2, \cdots, 4k - 1$, the necessary condition for the existence of solution of Problem \eqref{Pucci-Serrin} is that $\lambda$ should be larger than some positive constant number. Such dimensions are called critical dimensions. When $n = 3$ and $k = 1$, Brezis and Nirenberg \cite{BrezisNirenberg1} already found such a lower bound explicitly. Pucci and Serrin proved $n = 2k + 1$ is critical. When $k = 2$, the biharmonic version of Pucci-Serrin's conjecture has been proved by Edmunds, Fortunato and Jannelli \cite{EdmundsFortunatoJannelli1}. The cases $k = 3,4$ were due to Bernis and Grunau \cite{BernisGrunau1, Grunau1}.

\medskip

The main purpose of the present paper is to study the higher order Brezis-Nirenberg problem on hyperbolic spaces. To formulate our equation, we first briefly recall the proof of Brezis and Nirenberg to solve the Problem (\ref{Brezis-Nirenberg}) in \cite{BrezisNirenberg1}. If we define the functional

$$
  I[u] = \int_\Omega (|\nabla u|^2 - \lambda u^2 )dx
$$
and denote

$$
  S_\lambda = \inf_{u \in W^{1,2}_0(\Omega)\setminus \{ 0 \}} \frac{I_\lambda [u]}{\left( \int_\Omega |u|^q dx \right)^{2/q}}.
$$
where $q = \frac{2n}{n - 2}$, Brezis and Nirenberg discovered the following criterion:

\begin{theorem}[\cite{BrezisNirenberg1}]\label{Criterion}
  There exists at least one solution of Problem \eqref{Brezis-Nirenberg} when $\lambda$ satifsies

  $$
    S_\lambda < S_{n,1},
  $$
  where $S_{n,1}$ is the best constant of the first order Sobolev inequality, i.e.

  $$
    S_{n,1} = \inf_{u \in W^{1,2}_0(\Omega) \setminus \{ 0 \}}  \frac{ \int_\Omega |\nabla u|^2 dx }{\left( \int_\Omega |u|^q dx \right)^{2/q}}.
  $$
\end{theorem}
Sharp constants, together with extremal functions, of Sobolev inequalities have been well understood for decades. Among other results, Talenti \cite{Talenti1} and Aubin \cite{Aubin2} proved that for the sharp Sobolev inequality

$$
  S_{n,1} \left( \int_{\Bbb{R}^n} |u|^q dx \right)^{2/q} \leq \int_{\Bbb{R}^n} |\nabla u|^2 dx,
$$

$$
  S_{n,1} = \pi^{1/2} n^{-1/2} \left( \frac{1}{n - 2} \right)^{1/2} \left( \frac{\Gamma(1 + n/2) \Gamma(n)}{\Gamma(n/2) \Gamma(1 + n - n/2)} \right)^{1/n}.
$$
Moreover, the extremal function has the form (up to some translations and dilations)

$$
  u(x) = \frac{1}{( a + b |x|^2 )^{\frac{n-2}{2}}},
$$
where $a,b$ are positive constants. The knowledge of the best constant, together with the explicit form of extremal functions, plays an essential role in verifying the assumption of Theorem \ref{Criterion}. Precisely, in order to show $S_\lambda < S_{n,1}$, one needs to test a sequence in the form of the extremal function, with an appropriate smooth cut-off, to the quotient in $S_\lambda$ and perform an involved estimate. After the works of \cite{Talenti1} and \cite{Aubin2}, sharp Sobolev inequalities   together with their best constants and extremal functions have been studied in non-Euclidean spaces,  see e.g., \cite{Beckner1, DruetHebey1,Hebey1,  Liu1} and the references therein. These results hence enlighten the possibility to the study of Brezis-Nirenberg problems in non-Euclidean settings.
The sharp first order Sobolev inequalities of Aubin \cite{Aubin2} and Talenti \cite{Talenti1} have also been extended to higher order Sobolev inequalities in Euclidean spaces, see e,g.,  \cite{WangX},  \cite{CT},  \cite{Swanson}.
\medskip

Among numerous improvements of Sobolev inequalities, one has the following Hardy-Sobolev-Maz'ya's inequality on half spaces \cite{Mazya1}:

\begin{equation}\label{HSM}
  C \left( \int_{\Bbb{R}^n_+}  x_1^\gamma |u|^q dx \right)^{2/q} \leq \int_{\Bbb{R}^n_+} |\nabla u|^2 dx - \frac{1}{4}\int_{\Bbb{R}^n_+} \frac{|u|^2}{x_1^2} dx, \ u \in C_0^\infty (\Bbb{R}^n_+),
\end{equation}
where $\gamma = \frac{(n-2)q}{2} - n$ and $C=C(\gamma, n)$ is a positive constant. In particular, when $\gamma = 0$, this improves the classical Sobolev inequality due to the subtraction of an additional  Hardy term from the right hand side.

Recently, the second and third authors \cite{LuYang1} established the higher order version of inequality \eqref{HSM}. To be precise,
they proved the following
\begin{theorem}\label{Lu-Yang}
  Let $2 \leq k < n/2$ and $2 < q \leq \frac{2n}{n - 2k}$, there exits a positive $C$ such that for every $u \in C_0^\infty( \Bbb{R}^n_+ )$,

  $$
    C \left( \int_{\Bbb{R}^n_+} x_1^\gamma |u|^q dx	\right)^{2/q} \leq \int_{\Bbb{R}^n_+} |\nabla^k u|^2 dx - \prod_{i = 1}^k \frac{(2i - 1)^2}{4} \int_{\Bbb{R}^n_+} \frac{u^2}{x_1^{2k}} dx.
  $$
  where $\gamma = \frac{(n-2)q}{2} - n$. Moreover, the constant $\prod_{i = 1}^k \frac{(2i - 1)^2}{4}$ on the second term of the right hand side is best possible.
\end{theorem}
To prove the above inequality, new ideas have been developed: one needs to firstly transform the problem on the upper half space to a problem on the hyperbolic space $\Bbb{H}^n$ and make use of the Helgason-Fourier analysis on hyperbolic spaces  $\Bbb{H}^n$. It turns out that the inequality in Theorem \ref{Lu-Yang} is equivalent to the following inequality on $\Bbb{H}^n$:

\begin{equation}\label{Lu-Yang hyperbolic space}
  C \left( \int_{\Bbb{H}^n} |u|^q dV \right)^{2/q} \leq \int_{\Bbb{H}^n} P_k u \cdot u dV - \prod_{i = 1}^k \frac{(2i - 1)^2}{4} \int_{\Bbb{H}^n} u^2 dV,
\end{equation}
where $dV$ is the volume form on $\mathbb{H}^n$, $P_k = P_1 (P_1 + 2) \cdots (P_1 + k (k-1))$ is the GJMS operator (see \cite{GrahamJenneMasonSparling1,Juhl1}), $P_1 = - \triangle_{\Bbb{H}^n} - \frac{n(n-2)}{4}$ is the conformal Laplacian on $\Bbb{H}^n$ and $\triangle_{\Bbb{H}^n}$ is the Laplace-Beltrami operator. Such an equivalence has been observed when $k = 1$ in \cite{ManciniSandeep1}.

 \medskip

 Motivated by the work of \cite{LuYang1}, in the present paper, we will study the following $2k-$th order  Brezis-Nirenberg problem on the hyperbolic space $\mathbb{H}^n$:

\begin{equation}\label{Brezis-Nirenberg on bounded domain}
  \begin{cases}
    P_k u - \lambda u = |u|^{q - 2}u \ \ \ &\text{on } \Omega, \\
    \nabla^\alpha_{\Bbb{H}^n} u = 0 \ \ \ &\text{on } \partial \Omega, \\
    \alpha = 0, \cdots , k-1,
  \end{cases}
\end{equation}
where $q = \frac{2n}{n - 2k}$, $\Omega$ is a bounded domain in $\Bbb{H}^n$ with $C^1$ boundary and $\nabla_{\Bbb{H}^n}$ is the hyperbolic gradient (see the precise definition in Section 2). When $k = 1$, this problem has been studied in \cite{Benguria1,Stapelkamp1}. Among other results, it was shown that when $n \geq 4$, the Problem \ref{Brezis-Nirenberg on bounded domain} always has solution in $W^{1,2}_0(\Omega)$ when $0 < \lambda < \Lambda_1 (P_1 , \Omega)$, where $\Lambda_1 (P_1 , \Omega)$ is the first Dirichlet eigenvalue of $P_1$ on $\Omega$. In addition, when $n = 3$, $k = 1$ and $\Omega$ is a geodesic ball, there exists an constant $\tilde{\lambda} > 0$ such that Problem \ref{Brezis-Nirenberg on bounded domain} has a solution in $W^{1,2}_0(\Omega)$ if and only if $\tilde{\lambda} < \lambda < \Lambda_1 (P_1 , \Omega)$.

\medskip

To state our theorems, we recall the sharp $k$th-order Sobolev inequality in $\mathbb{R}^n$ for $n>2k$:
$$\int_{\mathbb{R}^n}|\nabla^k u|^2dx\geq S_{n,k} \big(\int_{\mathbb{R}^n}|u|^{\frac{2n}{n-2k}}dx\big)^{\frac{n-2k}{n}},$$
where the
 best Sobolev constant $S_{n, k}$  is given by
$$S_{n,k}=\frac{n(n-2)}{4}\left(\frac{n(n-2)}{4}-2\right) \left(\frac{n(n-2)}{4}-6 \right)\cdot\cdot\cdot \left(\frac{n(n-2)}{4}-k(k-1) \right)|\mathbb{S}^n|^{\frac{2k}{n}}.$$
 The equality holds if and only if $u$ takes the form as $$u(x)=\frac{1}{(1+|x-x_0|^2)^{\frac{n-2k}{2}}}$$
up to some translation and dilation. See e.g.,  Wang \cite{WangX}, Cotsiolis and  Tavoularis \cite{CT}, Swanson \cite{Swanson}.

\medskip

The first main theorem in this paper is the existence result of solutions to the Brezis-Nirenberg problem on the bounded domain in the hyperbolic space $\mathbb{H}^n$, namely Theorem
  \ref{Existence result at critical exponent on bounded domain} below. To show this, we will need
   a recent work from \cite{LuYang2},  where the second and third authors of the present paper proved the following for $k\ge 2$ (see also \cite{Hebey1} for the case $k = 1$):

\begin{theorem}\label{Lu-Yang theorem}Assume $2\le k<\frac{n}{2}$.
  If there exists a constant $\lambda \in \Bbb{R}$ such that for any $u \in C_0^\infty (\Bbb{H}^n)$,

  $$
    S_{n,k} \left( \int_{\Bbb{H}^n} |u|^q dV \right)^{2/q} \leq \int_{\Bbb{H}^n} P_k u \cdot u - \lambda u^2 dV,
  $$
where $S_{n,k}$ is the Sobolev best constant with $k$-th order derivative. Then the following are true:

 (i) When $n \geq 4k$, $\lambda \leq 0$ must hold; and

 (ii) When $2k+ 2 \leq n \leq 4k - 1$, $\lambda \leq \underline{\Lambda}$ must hold, where $\underline{\Lambda}$ can be taken as

 \begin{equation}\label{LambdaValue}
       \underline{\Lambda} = \frac{\Gamma(n/2) \Gamma(k) \sum_{j = 0}^{k-1} \frac{\Gamma(j + \frac{n-2k}{2})}{\Gamma(j + 1) \Gamma(\frac{n-2k}{2})}}{ 2^{\frac{n+2k}{2}} \Gamma(\frac{n-2k}{2}) \int_0^1 [2^{2k-1} - \sum_{j = 0}^{k-1} \frac{\Gamma(j + \frac{n-2k}{2})}{\Gamma(j + 1) \Gamma(\frac{n-2k}{2})} (1 - r^2)^j  ]^2 \frac{ r^{n-1} dr }{ (1 - r^2)^{2k} } }.
  \end{equation}
\end{theorem}

\begin{remark}
For the case $2k + 2 \leq n \leq 4k - 1$, Theorem \ref{Lu-Yang theorem} gives one value of $\underline{\Lambda}$ explicitly. The proof of this fact depends on a quite complicated calculation over a carefully designed testing function. This value of $\underline{\Lambda}$ can allow us to complete the proof of
the second part of Theorem \ref{Existence result at critical exponent on bounded domain} by using Theorem \ref{Lu-Yang theorem} directly.
However, the second statement of Theorem \ref{Existence result at critical exponent on bounded domain} can also be proved through a different argument and it works for  dimensions $2k + 1 \leq n \leq 4k - 1$. Moreover, in Section 7, we consider the special case $k = 2$ and $\Omega$ being a hyperbolic ball. We will obtain an explicit value of $\underline{\Lambda}$ which is smaller than the one given in \eqref{LambdaValue} and therefore sharpens the results of Theorem \ref{Lu-Yang theorem}.
\end{remark}

\medskip

Now we are ready to state
 the following existence of a solution to the higher order Brezis-Nirenberg problem  \eqref{Brezis-Nirenberg on bounded domain} on bounded domains in the hyperbolic space $\mathbb{H}^n$.

\begin{theorem}\label{Existence result at critical exponent on bounded domain}Assume $2\le k<\frac{n}{2}$.
Let $\Omega\subset \mathbb{H}^n$ be a bounded domain with $C^1$ boundary
and  $\Lambda_1(P_k, \Omega)$ be the first Dirichlet eigenvalue of $P_k$ on $\Omega$ defined by
$$
  \Lambda_1 (P_k, \Omega) = \inf_{u \in C_0^\infty (\Omega) \setminus \{ 0 \}} \frac{\int_\Omega P_ku\cdot u dV}{\int_\Omega |u|^2 dV}.
$$
Then
  \begin{itemize}
    \item For $n \geq 4k$ and $0 < \lambda < \Lambda_1 (P_k,\Omega)$, the $2k-$th order Brezis-Nirenberg problem  \eqref{Brezis-Nirenberg on bounded domain} has at least one nontrivial solution in $W^{k,2}_0 (\Omega)$ (the definition of Sobolev spaces $W^{k,2}_0 (\Omega)$ on $\Bbb{H}^n$ will be defined in Section 2).

    \item For $2k + 1 \leq n \leq 4k-1$, there exists a positive constant $\Lambda^*$ such that for $\Lambda^* < \lambda <  \Lambda_1(P_k, \Omega)$ the Brezis-Nirenberg problem \eqref{Brezis-Nirenberg on bounded domain} has at least one nontrivial solution in $W_0^{k,2}(\Omega)$.

    \end{itemize}
\end{theorem}

\begin{remark}
 One reason we require $\lambda < \Lambda_1(P_k , \Omega)$ is that the Hardy-Sobolev-Maz'ya inequality fails if $\lambda \geq \Lambda_1(P_k , \Omega)$ and our argument depends crucially  on such an inequality. Another reason is when $k = 1$ and $\lambda \geq \Lambda_1(P_1 , \Omega)$, one can easily prove the nonexistence of any positive solution to the second order Brezis-Nirenberg problem \eqref{Brezis-Nirenberg}, by testing a positive eigenfunction. However, such a method depends on the existence of a positive eigenfunction and it is not necessarily true for higher order equations. In fact, only partial results are known for higher order cases on balls centered at zero even in the Euclidean spaces (see \cite{GazzolaGrunauSweers1}).

\end{remark}

\begin{remark}
Moreover, for all the dimensions $2k+ 2 \leq n \leq 4k - 1$, if we denote $S_{n,k}$ as the sharp constant of the classical $k$-th order Sobolev inequality, $ \Lambda^*$ can be taken as
$$\Lambda^*=\min\left\{\Lambda_1 (P_k , \Omega) - |\text{Vol} (\Omega)|^{\frac{2}{q} - 1} S_{n,k},\,\, \underline{\Lambda}\right\}$$
where the constant $\underline{\Lambda}$ is given in Theorem \ref{Lu-Yang theorem}, and for $n=2k+1$, we can choose $\Lambda^*=\Lambda_1 (P_k , \Omega) - |\text{Vol} (\Omega)|^{\frac{2}{q} - 1} S_{n,k}$.
\end{remark}

\begin{remark}
The lower bound $\Lambda^*$ obtained in Theorem \ref{Existence result at critical exponent on bounded domain} above is most likely not optimal
for the existence of nontrivial solutions to
the $2k-$th order Brezis-Nirenberg problem  \eqref{Brezis-Nirenberg on bounded domain}. However, in the next theorem, we will show
an improved lower bound in the case of $k=2$ when the domain is a geodesic ball.
\end{remark}

\medskip

If we further consider the special case when $\Omega$ is a geodesic ball, note that for the ball model of $\Bbb{H}^n$, the Euclidean ball $B_R(0)$ (centered at origin with radius $R$) is also a geodesic ball on the hyperbolic space. Suppose the corresponding geodesic radius is $\overline{\rho}$, then they satisfy the relation $\overline{\rho} = \log \frac{1+R}{1 - R}$.  In lower dimensional cases when $k=2$ and $n=5, 6, 7$, we can improve the lower bound of $\lambda$ as follows:
\medskip

\begin{theorem}\label{Lower bound lambda}
  Assume $k = 2$ and  $4 < n < 8$ and $\Omega = B_R(0)$ is a Euclidean ball centered at the origin with radius $0<R<1$. If we define $\lambda^*$ as the first positive value such that

  $$
    \det
  \begin{pmatrix}
    P_{\lambda_1^*}^{\frac{2-n}{2}} (\overline{\rho}) & P_{\lambda_2^*}^{\frac{2-n}{2}} (\overline{\rho}) \\
    P_{\lambda_1^*}^{\frac{4-n}{2}} (\overline{\rho}) & P_{\lambda_2^*}^{\frac{4-n}{2}} (\overline{\rho})
  \end{pmatrix}
  = 0,
  $$
  where $\lambda_1^* = -1 + \sqrt{1 + \lambda^*}$, $\lambda_2^* = -1 - \sqrt{1 + \lambda^*}$, $\overline{\rho} = \log \frac{1+R}{1 - R}$ and $P_l^m$ is the Legendre function (see the precise definition in Section 7), then Problem \ref{Brezis-Nirenberg on bounded domain} has a  non-trivial solution when $\lambda^* < \lambda < \overline{\Lambda} (P_2 , \Omega)$.
\end{theorem}
\begin{remark}
 Recall the definition of $\Lambda^*$ in Remark 1.9, from the proof of Theorem \ref{Lower bound lambda}, we actually have $\lambda^* \leq \Lambda^*$. Moreover, we conjecture that this $\lambda^*$ is the sharp lower bound, i.e. the Problem \ref{Brezis-Nirenberg on bounded domain} has no solution when $\lambda \leq \lambda^*$. This conjecture remains to be further studied.
\end{remark}

\medskip

As for the nonexistence result, we further assume $\Omega$ is star-shaped, i.e. the inner product of the outer normal vector field and $\frac{\partial}{\partial \rho}$ is positive, where $\rho$ denotes the geodesic distance function on $\Bbb{H}^n$. Then the second main result of this paper reads as follows:

\begin{theorem}\label{Nonexistence} Assume $2\le k<\frac{n}{2}$.
Let $\Omega$ be a star-shaped domain in $\mathbb{H}^n$, the $2k-$th order Brezis-Nirenberg problem

$$
  \begin{cases}
  P_k u - \lambda u = |u|^{q - 2}u, \\
  \nabla_{\Bbb{H}^n}^\alpha u \vert_{\partial \Omega} = 0 \ \ \ \text{for } |\alpha| \leq k - 1
  \end{cases}
$$
has no nontrivial solution in $C^{2k} (\overline{\Omega})$ if
 $\lambda < 0$.

\end{theorem}
To prove Theorem \ref{Nonexistence}, we need the following identity (see \cite{Liu1, LuYang1}):

$$
   \left( \frac{1 - |x|^2}{2} \right)^{k + \frac{n}{2}} (- \triangle)^k \left( \left( \frac{1 - |x|^2}{2} \right)^{k - \frac{n}{2}} u \right) = P_k u,
$$
where $\Delta$ is the Laplace operator in the Euclidean space $\mathbb{R}^n$.

It implies the following equivalent higher order equation of \eqref{Brezis-Nirenberg on bounded domain}:

\begin{equation}\label{equivalent equation}
  (- \triangle)^k \tilde{u} = \tilde{u}^{q-1} + \lambda \left( \frac{2}{1 - |x|^2} \right)^{2k} \tilde{u},
\end{equation}
where $\tilde{u} = \left( \frac{1 - |x|^2}{2} \right)^{k - \frac{n}{2}} u$ and the equation is defined on a star-shaped domain inside of Euclidean ball. To this end, we further establish a Pohozaev type identity (see \cite{BrezisNirenberg1, GazzolaGrunauSweers1, Pohozaev1}), which implies Theorem \ref{Nonexistence}.

\medskip

We now consider the equation $P_k u - \lambda u = |u|^{q - 2} u $ on the entire hyperbolic space $\Bbb{H}^n$. When $k = 1$, this problem has been studied  in \cite{ManciniSandeep1}, where among other results,  they showed that the entire solution exists if and only if  $0< \lambda \leq  \frac{1}{4}$. However for higher order equations we observe several new phenomenons. More precisely, we prove the following existence result which is another main result of this paper.

\begin{theorem}\label{Existence result at critical exponent on whole space} Let $k\ge 2$.
  Consider the following equation on $\mathbb{H}^n$

  $$
    P_k u - \lambda u = |u|^{q - 2} u.
  $$
Then
  \begin{itemize}
    \item when $n \geq 4k$ and $0 < \lambda < \overline{\Lambda} \equiv \prod_{j = 1}^k \frac{(2j - 1)^2}{4} $,  there exists at least one nontrivial solution in $W^{k,2}(\Bbb{H}^n)$;

    \item when $2k + 2 \leq n \leq 4k-1$ and $\underline{\Lambda} < \lambda < \overline{\Lambda}$, there esxists at least one nontrivial solution in $W^{k,2}(\Bbb{H}^n)$, where $\underline{\Lambda}$ is defined as

    $$
      \underline{\Lambda} = \frac{\Gamma(n/2) \Gamma(k) \sum_{j = 0}^{k-1} \frac{\Gamma(j + \frac{n-2k}{2})}{\Gamma(j + 1) \Gamma(\frac{n-2k}{2})}}{ 2^{\frac{n+2k}{2}} \Gamma(\frac{n-2k}{2}) \int_0^1 [2^{2k-1} - \sum_{j = 0}^{k-1} \frac{\Gamma(j + \frac{n-2k}{2})}{\Gamma(j + 1) \Gamma(\frac{n-2k}{2})} (1 - r^2)^j  ]^2 \frac{ r^{n-1} dr }{ (1 - r^2)^{2k} } }
    $$
  \end{itemize}
\end{theorem}

Here are some remarks towards Theorem \ref{Existence result at critical exponent on whole space}. In the Euclidean case, Pucci and Serrin showed that $n = 2k + 1$ is a critical dimension, however such information is missing in Theorem \ref{Existence result at critical exponent on whole space}. When $n = 3, k = 1$, it was proved in \cite{ManciniSandeep1} that the
equation $P_k u - \lambda u = |u|^{q - 2} u $ on the whole space has no solution for any $\lambda$.  In the proof of Theorem \ref{Existence result at critical exponent on bounded domain}, we will establish a criterion which asserts that the sharp constant of Hardy-Sobolev-Maz'ya inequality being strictly less than the Sobolev best constant will imply the existence of solution to the Brezis-Nirenberg problem. Combing this with Theorem \ref{Lu-Yang theorem} will partially explain the missing dimension and as for higher orders, we tend to believe that whatever value $\lambda$ takes, there is no solution when $n = 2k+1$.

\medskip

The symmetry of solutions plays an essential role in the study of Brezis-Nirenberg problem, especially in the lower dimension case. From \cite{BrezisNirenberg1}, by proving the symmetry of solutions to Problem \ref{Brezis-Nirenberg}, the study of the solution gap problem is reduced to the study of the ODE generated from the original equation (see also \cite{Benguria1, ManciniSandeep1} for non-Euclidean cases). In the present paper, we will prove the following symmetry result for the higher order equations:

\begin{theorem}\label{symmetry result} Let $k\ge 2$.
  If $u \in W_0^{k,2}(\Bbb{H}^n)$ is a positive weak solution of the equation

  \begin{equation}\label{differential equation}
    P_k u - \lambda u = |u|^{q - 2}u,
  \end{equation}
   then there exists a point $P \in \Bbb{H}^n$ such that $u$ is constant on the geodesic spheres centered at $P$. Moreover, $u$ is nonincreasing.
\end{theorem}

Before we continue, we will have some discussions on symmetry of solutions to classical elliptic differential equations in the Euclidean spaces.
The symmetry of solutions for semilinear elliptic problem on $\Omega\subset \mathbb{R}^n$:

$$
  \begin{cases}
    - \triangle u = f(u) \ \text{on} \ \Omega ,\\
    u = 0 \ \text{on} \ \partial \Omega
  \end{cases}
$$
was studied by Gidas, Ni and Nirenberg \cite{GidasNiNirenberg1}, using the moving plane method. Such a method was initiated by Alexandroff, and was further developed by Serrin \cite{Serrin1}, Gidas, Ni and Nirenberg \cite{GidasNiNirenberg1}, Caffarelli, Gidas and Spruck \cite{CaffarelliGidasSpruck1}. In the work of Almeida, Damascelli and Ge \cite{AlmeidaDamascelliGe1} (see also \cite{AlmeidaGe1}), such method was extended to manifolds with group symmetry (this terminology will be rigorously defined in the later section). One specific example of such a family of manifolds is the hyperbolic space $\Bbb{H}^n$.  However, it is impossible to directly apply the classical moving plane method to higher order equations, due to the fact that such method depends heavily on the maximum principle, which are not available for higher order cases. If one further assumes that the nonlinearity $f(t) : [0 , \infty) \to \Bbb{R}$ is continuous, nondecreasing and $f(0) \geq 0 $, $\Omega$ is a Euclidean ball, Berchio, Gazzola and Weth \cite{BerchioGazzolaWeth1} showed that a modified moving plane argument would still give symmetry result for polyharmonic Dirichlet problems. Unfortunately, such a method fails on $\Bbb{H}^n$ since the volume on $\Bbb{H}^n$ has an exponential growth, which makes the nonlinearity term even worse. On the other hand, Chen, Li and Ou \cite{ChenLiOu1} developed an integral equation version of the moving plane argument.   With the help of the Green's function estimates obtained in \cite{LuYang1} and the Helgason-Fourier transform on hyperbolic spaces, we will prove that the higher order Brezis-Nirenberg problem on the entire space $\mathbb{H}^n$ is equivalent to the integral equations and then establish the symmetry of the solutions by developing a moving plan argument in integral form in the hyperbolic spaces in the spirit of Chen, Li and Ou \cite{ChenLiOu1}.

\medskip
Therefore, to prove Theorem \ref{symmetry result}, we will show that the higher order differential equation (\ref{differential equation}) is equivalent to an integral equation using the Helgason-Fourier analysis on the hyperbolic spaces.

\begin{theorem}\label{Integral equation theorem}  If $u \in W_0^{k,2}(\mathbb{B}^n)$ is a positive weak solution of the higher order differential equation
(\ref{differential equation}), then $u$ must satisfy the following integral equation for the Green function $G(x, y)$ of the operator $P_k - \lambda$
on the hyperbolic ball $\mathbb{B}^n$:
  \begin{equation}\label{Integral equation1}
     u(x) = \int_{\Bbb{B}^n} G(x,y) u^{q-1}(y) dV_y.
  \end{equation}
Moreover, for any  fixed $y\in \mathbb{B}^n$,  $G(x,y)$ is a positive radially decreasing function with respect to the geodesic distance $\rho = d(x,y)$.

\end{theorem}

We end our introduction with the following remarks. There are substantial differences in our approach to study the Brezis-Nirenberg problem on the hyperbolic spaces  from those in the Euclidean spaces. First of all, we rely on the higher order Hardy-Sobolev-Maz'ya inequalities on hyperbolic spaces established by the second and third authors \cite{LuYang1} and \cite{LuYang2}.
To the best of our knowledge, so far very little is known towards the higher order Brezis-Nirenberg problem  in the non-Euclidean settings (e.g. Problem \ref{Brezis-Nirenberg on bounded domain}). This is partially because such a problem is motivated by the higher order Hardy-Sobolev-Maz'ya inequality, which is merely proved very recently by the second and third authors. There are still many open questions and we want to make some remarks here. As mentioned previously, in the Euclidean space $\mathbb{R}^3$, the Brezis-Nirenberg Problem \ref{Brezis-Nirenberg} on a ball has solution if and only if $\frac{1}{4} \Lambda_1(- \triangle, \Omega) < \lambda < \Lambda_1(- \triangle , \Omega)$. Such a \textbf{solution gap} problem has also been  solved for the second order Brezis-Nirenberg problem on geodesic balls in $\Bbb{H}^n$ by S. Benguria in \cite{Benguria1}. For higher order equations on $\Bbb{H}^n$, even when in the lower order cases  $k = 2,3$, the solution gap issue remains to be further studied. Though we have found a lower bound $\lambda^*$ for $k= 2$ such that the Brezis-Nirenberg problem on any geodesic ball $\Omega=B_R(O)$ in $\mathbb{H}^n$ for $n=5, 6, 7$  is solvable for $\lambda^*<\lambda <\Lambda_1(P_k, \Omega)$, we have not been able to show  yet if our lower bound $\lambda^*$ is sharp. (see Theorem \ref{Lower bound lambda}).
Second,  the Helgason-Fourier analysis and Green's function estimates on the hyperbolic spaces play an important role in our proofs. These are reflected in the proofs of  Theorem \ref{Existence result at critical exponent on bounded domain}, Theorem \ref{Lower bound lambda} and Theorem \ref{Integral equation theorem}.
Third, we prove the symmetry of solutions by converting the higher order differential equations to the integral equations using the Helgason-Fourier analysis on the hyperbolic spaces so that we can avoid the maximum principle which is not available for the higher order differential equations. Our proof of symmetry is also substantially different from that in the first order equations on the hyperbolic spaces \cite{ManciniSandeep1} where the equations are converted
into the corresponding ones on the Euclidean space and that method does not apply to the higher order case of our study.

\medskip

The present paper will be organized as follows: in Section 2, we will give some preliminaries and useful known facts; Section 3 will be devoted to the proof of Theorem \ref{Existence result at critical exponent on bounded domain}; Section 4 will give the proof of Theorem \ref{Existence result at critical exponent on whole space}; Section 5 will focus on Theorem \ref{Nonexistence}; Section 6 deals with the proof of Theorem \ref{Integral equation theorem} and Theorem \ref{symmetry result}; Section 7 will focus on sharpening  Theorem \ref{Existence result at critical exponent on bounded domain} by improving its bound of $\lambda$ and then prove Theorem \ref{Lower bound lambda}.

\section{Notations and preliminaries}
We begin by quoting some preliminaries which will be needed in the sequel and refer to \cite{Ahlfors1,AlmeidaDamascelliGe1,GelfandGindikin1,Helgason1,Helgason2,Hua1,LiuPeng1} for more information about this subject.

\subsection{The half space model of $\Bbb{H}^n$}

It is given by $\Bbb{R}^{n-1} \times \Bbb{R}_+ = \{ (x_1 , \cdots , x_n) : x_1 > 0 \}$ equipped with the Riemannian metric $ds^2 = \frac{dx_1^2 + \cdots + dx_n^2}{x_1^2}$. The volume form is $dV = \frac{dx}{x_1^n}$, where $dx$ is the Lebesgue measure on $\Bbb{R}^n$. The hyperbolic gradient is $\nabla_\Bbb{H} = x_1 \nabla$ and the Laplace-Beltrami operator on $\Bbb{H}^n$ is given by

$$
  \triangle_\Bbb{H} = x_1^2 \triangle - (n-2) x_1 \frac{\partial}{\partial x_1},
$$
where $\triangle$ is the usual Laplacian on $\Bbb{R}^n$.

\subsection{The ball model of $\Bbb{H}^n$}

It is given by the unit ball $\Bbb{B}^n$ equipped with the usual Poincar\'e metric

$$
  ds^2 = \frac{4 (dx_1^2 + \cdots + dx_n^2)}{(1 - |x|^2)^2}.
$$
The volume form is $dV = \left( \frac{2}{1 - |x|^2} \right)^n dx$. The hyperbolic gradient is $\nabla_{\Bbb{H}^n} = \frac{1 - |x|^2}{2} \nabla$ and the Laplace-Beltrami operator is given by

$$
  \triangle_{\Bbb{H}^n} = \frac{1 - |x|^2}{4} \left( (1 - |x|^2) \triangle + 2(n - 2) \sum_{i = 1}^n x_i \frac{\partial}{\partial x_1} \right).
$$
Furthermore, the half space model and ball model are equivalent.

\subsection{Sobolev spaces on hyperbolic space $\Bbb{H}^n$}

We will define Sobolev space $W^{k,2}$ on the ball model. For any open $\Omega \subset \Bbb{H}^n$ and $u \in C^\infty (\Omega)$, set $|\nabla_{\Bbb{H}^n} u|^2 = \left( \frac{1 - |x|^2}{2} \right)^2 |\nabla u|$. Define the $W^{k,2}$-norm of $u$ as

$$
  ||u||_{W^{k,2} (\Omega)} = \sum_{0 \leq j \leq k, j \textit{ is even}} \int_\Omega |(- \triangle_{\Bbb{H}^n})^{j/2} u|^2 dV + \sum_{0 \leq j \leq k, j \textit{ is odd}} \int_\Omega |\nabla_{\Bbb{H}^n} (- \triangle_{\Bbb{H}^n})^{\frac{j-1}{2}} u|^2 dV.
$$
Then the Sobolev space $W^{k,2}(\Omega)$ is the closure of $C^\infty (\Omega)$ with respect to $||\cdot||_{W^{k,2}(\Omega)}$. In addition, the Sobolev space $W^{k,2}_0(\Omega)$ is the closure of $C^\infty_0 (\Omega)$ with respect to $||\cdot||_{W^{k,2}(\Omega)}$.

\subsection{Foliations of hyperbolic spaces $\Bbb{H}^n$}\label{section-foliation}

In this subsection, we introduce the foliation of $\Bbb{H}^n$ which is needed in the study of the hyperbolic symmetry of the solution (see \cite{AlmeidaDamascelliGe1,AlmeidaGe1} for more detail). Let $\Bbb{R}^{n,1} = ( \Bbb{R}^{n+1} , g )$ be the Minkowski space, where the metric $ds^2 = -dx_0^2 + dx_1^2 + \cdots + dx_n^2$. The hyperboloid model of hyperbolic space $\Bbb{H}^n$ is the submanifold $\{ x \in \Bbb{R}^{n,1}: g(x,x) = -1, x_0 > 0 \}$. A particular directional foliation can be obtained by choosing any direction in the $x_1, x_2, \cdots , x_n$ plane. Without loss of generality, we consider $x_1$ direction. Denote $\Bbb{R}^{n,1} = \Bbb{R}^{1,1} \times \Bbb{R}^{n-1}$, where $(x_1, x_0) \in \Bbb{R}^{1,1}$. Define $A_t = \tilde{A}_t \otimes Id_{\Bbb{R}^{n-1}}$, where $\tilde{A}_t$ is the hyperbolic rotation on $\Bbb{R}^{1,1}$:

$$
  \tilde{A}_t = \left(\begin{array}{cc}\cosh t & \sinh t \\\sinh t & \cosh t \end{array}\right).
$$
The reflection is defined by $I (x_0, x_1, x_2,  \cdots, x_n) = (x_0, -x_1, x_2, \cdots, x_n)$. Let $U = \Bbb{H}^n \cap \{ x_1 = 0 \}$ and $U_t = A_t (U)$, then $\Bbb{H}^n$ is foliated by $U_t$, i.e. $\Bbb{H}^n = \cup_{t \in \Bbb{R}} U_t$. Moreover, if we define $I_t = A_t \circ I \circ A_{-t}$, then it is easy to verify that $I_t (U_t) = U_t$.

\subsection{Helgason-Fourier transform on the hyperbolic space $\Bbb{H}^n$}
We recall here the Helgason-Fourier analysis on hyperbolic spaces. We refer the reader to
\cite{Helgason1} and \cite{Helgason2}.

Set

$$
  e_{\lambda, \zeta} (x) = \left( \frac{\sqrt{1 - |x|^2}}{|x - \zeta|} \right)^{n - 1 + i \lambda}, \ x \in \Bbb{B}^n , \ \lambda \in \Bbb{R}, \  \zeta \in \Bbb{S}^{n-1}.
$$
The Fourier transform of a function $u$ on $\Bbb{H}^n$ (ball model) can be defined as

$$
  \hat{u}(\lambda , \zeta) = \int_{\Bbb{B}^n} u(x) e_{-\lambda , \zeta}(x) dV,
$$
provided this integral exists. The following inversion formula holds for $u \in C_0^\infty (\Bbb{B}^n)$:

$$
  u(x) = D_n \int_{- \infty}^{+ \infty} \int_{\Bbb{S}^n} \hat{u} (\lambda, \zeta) e_{\lambda, \zeta} (x) |\mathfrak{c}(\lambda)|^{-2} d\lambda d\sigma(\zeta),
$$
where $D_n = \frac{1}{2^{3 - n} \pi |\Bbb{S}^{n-1}|}$ and $\mathfrak{c} (\lambda)$ is the Harish-Chandra $\mathfrak{c}$-function given by

$$
  \mathfrak{c} (\lambda) = \frac{2^{n - 1 - i \lambda} \Gamma(n/2) \Gamma(i \lambda)}{\Gamma(\frac{n - 1 + i \lambda}{2}) \Gamma(\frac{1 + i \lambda}{2})}.
$$
Similarly, there holds the Plancherel formula:

$$
  \int_{\Bbb{H}^n} |u (x)|^2 dV = D_n \int_{-\infty}^{+\infty} \int_{\Bbb{S}^{n-1}} |\hat{u} (\lambda, \zeta)|^2  |\mathfrak{c} (\lambda)|^{-2} d\lambda d\sigma(\zeta).
$$
Since $e_{\lambda, \zeta}(x)$ is an eigenfunction of $-\triangle_{\Bbb{H}^n}$ with eigenvalue $\frac{(n-1)^2 + \lambda^2}{4}$, it is easy to check that for $u \in C_0^\infty (\Bbb{H}^n)$,

$$
  \widehat{\triangle_{\Bbb{H}^n} u} (\lambda,\zeta) = - \frac{(n-1)^2 + \lambda^2}{4} \hat{u}(\lambda,\zeta).
$$
Therefore, in analogy with the Euclidean setting, we define the fractional Laplacian on hyperbolic spaces as following:

$$
  \widehat{(- \triangle_{\Bbb{H}^n})^\gamma u} (\lambda,\zeta) = \left( \frac{(n-1)^2 + \lambda^2}{4}  \right)^\gamma \hat{u} (\lambda,\zeta).
$$

\subsection{Some known facts}

The spectral gap of $- \triangle_{\Bbb{H}^n}$ on $\L^2 (\Bbb{H}^n)$ is $\frac{(n-1)^2}{4}$, i.e.

$$
  \frac{(n-1)^2}{4} \int_{\Bbb{H}^n} u^2 dV \leq \int_{\Bbb{H}^n} |\nabla_{\Bbb{H}^n} u|^2 dV, \ u \in C_0^\infty (\Bbb{H}^n).
$$
In \cite{Liu1}, Liu proved the following sharp Sobolev inequality:

\begin{equation}\label{Liu}
  S_{n,k} \left( \int_{\Bbb{H}^n} |u|^{\frac{2n}{n - 2k}} dV \right)^{\frac{n - 2k}{n}} \leq \int_{\Bbb{H}^n} P_k u \cdot u dV, \ u \in C_0^\infty (\Bbb{H}^n) \ \text{and} \ 1 \leq k < n/2.
\end{equation}
It is easy to see that \eqref{Lu-Yang hyperbolic space} improves \eqref{Liu} and Theorem \ref{Lu-Yang theorem} shows the criterion to achieve the sharp constant. The proof of  \eqref{Lu-Yang hyperbolic space}  depends on the Green's function estimate associated to the operator $P_k ^{-1}$. To be precise, with the help of the heat kernel on hyperbolic spaces, one has the following formula:

$$
  P_1^{-1} = \left(- \triangle_{\Bbb{H}} - \frac{n (n-2)}{4} \right)^{-1} = \frac{1}{n (n-2) \alpha(n)} \left( \frac{1}{(2\sinh \rho/2)^{n-2}} - \frac{1}{(2 \cosh \rho/2)^{n-2}} \right),
$$
where $\alpha(n) = \frac{ \pi^{n/2} }{ \Gamma(n/2 +1) } $ and $\rho$ denotes the geodesic distance. Moreover, in \cite{LuYang1}, the second and third authors proved

\begin{equation}\label{Green function estimate}
  G(x,y) \lesssim \left( \frac{1}{\sinh \frac{d(x,y)}{2}} \right)^{n - 2k},
\end{equation}
where $d(x,y)$ is the geodesic distance for $x,y \in \Bbb{H}^n$. This estimate enables us to apply the Hardy-Littlewood-Sobolev inequality on hyperbolic spaces (see Beckner \cite{Beckner2}, also \cite{LuYang1} for a proof):

\begin{theorem}\label{Hardy-Littlewood-Sobolev on H^n}
  Let $0 < \lambda < n$ and $p = \frac{2n}{2n - \lambda}$. Then for $u,v \in L^p (\Bbb{H}^n)$,

  $$
    \vert \int_{\Bbb{H}^n \times \Bbb{H}^n} \frac{u(x) v (y)}{( 2 \sinh \frac{d(x,y)}{2} )^\lambda} dV_x dV_y \vert \leq C_{n,\lambda} ||u||_p ||v||_p,
  $$
  where

  $$
     C_{n,\lambda} = \pi^{\lambda/2} \frac{\Gamma(n/2 - \lambda/2)}{\Gamma(n - \lambda/2)} \left( \frac{\Gamma(n/2)}{\Gamma(n)} \right)^{-1 + \lambda/n}
  $$
  is the best Hardy-Littlewood-Sobolev constant on $\Bbb{R}^n$. Furthermore, the constant $C_{n,\lambda}$ is sharp and there is no nonzero extremal function.
\end{theorem}

\section{Existence of solutions on the bounded domain: Proof of Theorem \ref{Existence result at critical exponent on bounded domain}}

The main purpose of this section is to establish the existence of solutions to the Brezis-Nirenberg problem on bounded domains in hyperbolic spaces, namely Theorem \ref{Existence result at critical exponent on bounded domain}. It is interesting to note that the Helgason-Fourier analysis on the hyperbolic spaces plays an important role in our approach given in this section.

\textbf{Proof of Theorem \ref{Existence result at critical exponent on bounded domain}}: Due to the density of $C_0^\infty (\Omega)$ in $W_0^{k,2}(\Omega)$, without loss of generality, we consider a minimizing sequence $\{ u_m \}$ in $C_0^\infty (\Omega)$ satisfying $||u_m||_{L^q (\Omega)} = 1$ and

$$
I_{\lambda,k} [u_m] = \int_{\Bbb{H}^n} P_k u_m \cdot u_m dV- \lambda u_m^2 dV \to S_{\lambda,k} \ \   \text{as} \ \ m \to \infty.
$$
We are to show that $\{ u_m \}$ is bounded in $W^{k,2}_0 (\Omega)$ when $\lambda<\Lambda_1 (P_k, \Omega)$. Recall that
$\Lambda_1 (P_k, \Omega)$ is the first eigenvalue of $P_k$ on $\Omega$.  Indeed, we have

$$
  (\Lambda_1 (P_k, \Omega) - \lambda) \int_\Omega u_m^2 dV \leq \int_\Omega P_k u_m \cdot u_m - \lambda u_m^2 dV \to S_{\lambda,k},
$$
hence $||u_m||_{L^2(\Omega)}$ is uniformly bounded. On the other hand, for any $j$-th order derivative, where $1 \leq j \leq k$, when $j$ is even,
by the Plancherel formula on the hyperbolic space $\mathbb{H}^n$ we have
\begin{align*}
  ||\triangle_{\Bbb{H}^n}^{j/2} u_m||_{L^2 (\Omega)} ^2
    &=  \int_\Omega \triangle_{\Bbb{H}^n}^j u_m \cdot u_m dV \\
    &= D_n \int_{- \infty}^{\infty} \int_{\Bbb{S}^{n-1}} \left( \frac{(n - 1)^2 + r^2}{4} \right)^j | \hat{u}_m(r , \zeta) |^2 | \mathfrak{c}(r) |^{-2} dr d\sigma(\zeta). \\
\end{align*}
We claim that when $\lambda<\Lambda_1 (P_k, \Omega)$, we have

$$
  \left( \frac{(n - 1)^2 + r^2}{4} \right)^j \lesssim \prod_{\ell = 1}^k \left( \frac{ (2\ell - 1)^2 + r^2 }{4} \right) - \lambda.
$$
Indeed, let $0 < \lambda <\Lambda_1 (P_k, \Omega)$ be fixed. When $r$ is large, since $j \leq k$, our claim obviously holds. When $r$ is close to zero, both left hand side and right hand side approach to fixed constants. Thus our claim is proved. This implies

\begin{align*}
   ||\triangle_{\Bbb{H}^n}^{j/2} u_m||_{L^2 (\Omega)} ^2
     &\lesssim  D_n \int_{- \infty}^{\infty} \int_{\Bbb{S}^{n-1}} \left( \prod_{\ell = 1}^k \left( \frac{ (2\ell - 1)^2 + r^2 }{4} \right) - \lambda  \right) | \hat{u}_m(r , \zeta) |^2 | \mathfrak{c}(r) |^{-2} dr d\sigma(\zeta), \\
     &= \int_\Omega P_k u_m \cdot u_m - \lambda u_m^2 dV ,\\
     &\to S_{\lambda, k}.
\end{align*}
When $j$ is odd, we can prove the uniform boundedness of $||\nabla_{\Bbb{H}} \triangle_{\Bbb{H}}^{\frac{j - 1}{2}} u_m||_{L^2 (\Omega)}$ in a similar manner. Therefore $\{ u_m \}$ is bounded in $W_0^{k,2}(\Omega)$ and hence $\{ u_m \}$ converges weakly to some $u \in W^{k,2}_0 (\Omega)$, and strongly in $L^2 (\Omega)$ (the compactness embedding comes from $C^1$ boundary condition, see T. Aubin \cite{Aubin1}). We do not always have the compactness with the critical exponent, but we are to prove that the compactness is achieved for some specific range of $\lambda$.
More precisely, note that

\begin{align}
  S_{\lambda,k} &= I_{\lambda,k} [u_m] + o(1) \nonumber\\
                    &= \int_{\Bbb{H}^n} \sum_{j = 1}^k c_j |\nabla_{\Bbb{H}^n}^j (u_m - u)|^2 dV + I_{\lambda,k} [u] + o(1) \nonumber\\
                    &\geq  S_{n,k} ||u_m - u||_q^2 + S_{\lambda,k} ||u||_q^2 + o(1) \\
                    &= (S_{n,k} - S_{\lambda,k}) ||u_m - u||_q^2  + S_{\lambda,k} ( ||u||_q^2 + ||u_m - u||_q^2 )  + o(1)\nonumber \\
                    &\geq (S_{n,k} - S_{\lambda,k}) ||u_m - u||_q^2 + S_{\lambda,k} (1 + o(1)) + o(1)\nonumber
\end{align}
where $S_{n,k}$ is the classical Sobolev best constant. Once we verify each line above, then as long as $S_{\lambda,k} < S_{n,k}$, we have the strong convergence of $\{ u_m \}$ in $L^q (\Omega)$.

From the first line to the second line in the above string of inequalities, it is due to the fact that

$$
  \int_{\Bbb{H}^n }P_k u \cdot u dV = \int_{\Bbb{H}^n} \sum_{j = 1}^k c_j |\nabla_{\Bbb{H}^n}^j u|^2 dV \ \  \text{for some constants } c_j
$$
and an application of the $W_0^{k,2}$ weak convergence:

$$
  \int_{\Bbb{H}^n} |\nabla^j_{\Bbb{H}^n} (u_m - u)|^2 dV = \int_{\Bbb{H}^n} |\nabla^j_{\Bbb{H}^n} u_m|^2 dV - \int_{\Bbb{H}^n} |\nabla^j_{\Bbb{H}^n} u|^2 dV + o(1).
$$
From the second line to the third line in the above string of inequalities, we use the Poincar\'e-Sobolev inequality \eqref{Liu}. The fourth line to the fifth line comes from the Brezis-Lieb lemma (see \cite{BrezisLieb1}) and the fact that $q \geq 2$:

$$
  ||u_m||_q^2 \leq ||u||_q^2 + ||u_m - u||_q^2 + o(1).
$$
\medskip

Hence, in order to establish the existence of a solution of the Brezis-Nirenberg problem on a bounded domain $\Omega$, the only question is: when would $S_\lambda < S_{n,k}$ hold? However this is completely characterized by Theorem \ref{Lu-Yang theorem}. In fact, Theorem \ref{Lu-Yang theorem} immediately tells us that we have  $S_\lambda < S_{n,k}$ for $\lambda<0$ when $n \geq 4k$, and $\lambda>\underline{\Lambda}$ for $2k+ 2 \leq n \leq 4k - 1$, where $\underline{\Lambda}$ is  defined as in Theorem \ref{Lu-Yang theorem}. Therefore, we have shown the existence
of a solution to the Brezis-Nirenberg problem on the bounded domain $\Omega$.

\medskip
 Nevertheless, we will further show below that we can actually solve the Brezis-Nirenberg problem for all dimensions  $2k+ 1 \leq n \leq 4k - 1$, pushing down the dimension to $n=2k+1$. To this end, we first denote $\phi$ as the eigenfunction of $P_k$ associated to the first eigenvalue $\Lambda_1 (P_k, \Omega)$. The existence of such function follows directly from a variational argument to

$$
  \Lambda_1 (P_k, \Omega) = \inf_{u \in W^{k,2}_0 (\Omega)} \frac{\int_{\Bbb{H}^n} \sum_{j = 1}^k c_j |\nabla_{\Bbb{H}^n}^j u |^2 dV}{||u||_2^2}.
$$
Moreover, we can assume $||\phi||_q = 1$ and we have

\begin{align*}
  S_{\lambda,k} &\leq I_{\lambda,k} [\phi] \\
                         &= (\Lambda_1(P_k, \Omega) - \lambda) \int_{\Bbb{H}^n} |\phi|^2 dV \\
                         &< (\Lambda_1(P_k, \Omega) - \underline{\Lambda}) \left( \int_{\Bbb{H}^n} |\phi|^q dV \right)^{2/q} |\text{Vol}(\Omega)|^{1 - \frac{2}{q}}\\=& (\Lambda_1(P_k, \Omega) - \underline{\Lambda})  |\text{Vol}(\Omega)|^{1 - \frac{2}{q}}.
\end{align*}
Then we just need to set $\lambda> \Lambda_1 (P_k , \Omega) - |\text{Vol} (\Omega)|^{\frac{2}{q} - 1} S_{n,k}$ and the above inequalities will give us $S_{\lambda,k} < S_{n,k}$.

In conclusion, for all the dimensions $2k+ 2 \leq n \leq 4k - 1$, we can choose
$$\Lambda^*=\min\{\Lambda_1 (P_k , \Omega) - |\text{Vol} (\Omega)|^{\frac{2}{q} - 1} S_{n,k},\,\, \underline{\Lambda}\}$$
and for $n=2k+1$, we can choose $\Lambda^*=\underline{\Lambda}$. Therefore, for $\lambda>\Lambda^*$,
the higher order Brezis-Nirenberg problem  \eqref{Brezis-Nirenberg on bounded domain} on bounded domains in the hyperbolic space $\mathbb{H}^n$
 has a solution.

\section{Existence of solutions on the whole space: Proof of Theorem \ref{Existence result at critical exponent on whole space}}

In this section, we will show the existence of the Brezis-Nirenberg problem on the entire hyperbolic space $\mathbb{H}^n$, namely Theorem \ref{Existence result at critical exponent on whole space}.

\textbf{Proof of Theorem \ref{Existence result at critical exponent on whole space}}:

Consider the Nehari manifold

$$
  \mathcal{S} = \{ u \in C_0^\infty(\mathbb{H}^n) : \ I_{\lambda,k} [u] < \infty, u \neq 0, I_{\lambda,k} [u] = \int_{\Bbb{H}^n} |u|^q dV \}.
$$
Notice that $$S_{\lambda,k} = \inf_{u \in \mathcal{S} } I_{\lambda,k} [u] ^{\frac{q-2}{q}} = \inf_{u \in \mathcal{S} } ||u||_q^{q-2}.$$
 Let $u_j$ be a minimizing sequence in $\mathcal{S}$, we first perform a proper dilation and translation. Given $(x_1,y) \in \Bbb{H}^n$ with the upper half space model ($x_1 \in \Bbb{R}_+, y \in \Bbb{R}^{n-1}$). For $y_0 \in \Bbb{R}^{n-1}$, let $B (y_0,R) = \{ (x_1,y) \in \Bbb{H}^n : x_1^2 + |y - y_0|^2 < R^2 \} $. Denote the concentration function $Q_j(R) = \sup_{y_0 \in \Bbb{R}^{n-1}} \int_{B (y_0, R)} |u_j|^q dV$, then for any $0 < \delta < S_{\lambda,k}^{\frac{q}{q-2}}$, one can always find $y_j \in \Bbb{R}^{n-1}$ and $R_j > 0$ such that

$$
  \delta = \int_{B (y_j, R_j)} |u_j|^q dV = \sup_{y_0 \in \Bbb{R}^{n-1}} \int_{B (y_0, R_j)} |u_j|^q dV.
$$

\noindent Let $v_j (x_1,y) = u_j ( (0,y_j) + R_j (x_1,y))$. Then $\{ v_j \} \subset \mathcal{S}$ and is still minimizing. To see this, one needs to check that

$$
  ||u_j||_2 = ||v_j||_2, ||u_j||_q = ||v_j||_q, I_{\lambda,k} [u_j] = I_{\lambda,k} [v_j] \to S_{\lambda,k}^\frac{q}{q-2}.
$$
\noindent To verify the last equality, notice that

$$
  \triangle_{\Bbb{H}^n} = x_1^2 \triangle - (n-2)x_1 \partial_{x_1},
$$

\noindent and hence a  change of variables gives

$$
  \int_{\Bbb{H}^n} (-\triangle_{\Bbb{H}^n})^{\ell} u_j dV = \int_{\Bbb{H}^n} (-\triangle_{\Bbb{H}^n})^{\ell} v_j dV, \text{ for any } 1 \leq \ell \leq k.
$$

\noindent Therefore, from the definition of GJMS operator $P_k$, we see that $\int_{\Bbb{H}^n} P_k u_j \cdot u_j dV = \int_{\Bbb{H}^n} P_k v_j \cdot v_j dV$. Moreover,

$$
 \delta = \int_{B(0,1)} |v_j|^q dV = \sup_{y_0 \in \Bbb{R}^{n-1}} \int_{B(y_j, 1)} |v_j|^q dV.
$$

\noindent By the Ekeland principle, we have

\begin{equation}\label{P-S condition}
   < v_j , u>_\lambda = \int_{\Bbb{H}^n} P_k v_j \cdot u - \lambda v_j u  dV = \int_{\Bbb{H}^n} |v_j|^{q-2} v_j u dV + o(1)
\end{equation}

\noindent for any $u$ in a bounded set of $W_0^{k,2}(\Bbb{H}^n)$. Following a similar argument as in the proof of Theorem \ref{Existence result at critical exponent on bounded domain}, we assert that $\{ v_j \}$ is a bounded sequence in $W_0^{k,2}(\Bbb{H}^n)$. Hence $v_j \rightharpoonup v \in W_0^{k,2}(\Bbb{H}^n)$. Moreover, $v_j \to v$ strongly in $L^p_{\text{loc}} (\Bbb{H}^n)$ and hence almost everywhere. Choosing $u = v$, we have $I_{\lambda,k} [v] = \int_{\Bbb{H}^n} |v|^q dV$. Therefore, it remains to show that $v \neq 0$. Assume by contradiction that $v = 0$.  Let $\phi \in C_0^\infty (B( (1,0) , R ))$ satisfy $\phi = 1$ on $B( (1,0) , R/2 )$ and $ 0 \leq \phi \leq 1$. Using $\phi^2 v_j$ as test function in \eqref{P-S condition}, we get

$$
   < v_j , \phi^2 v_j >_\lambda = \int_{\Bbb{H}^n} |v_j|^{q-2} (\phi v_j)^2 dV + o(1).
$$

\noindent On the other hand,

\begin{align*}
  &\int_{\Bbb{H}^n} P_k (\phi v_j) (\phi v_j) dV - \lambda (\phi v_j)^2 dV \\
  =& \int_{\Bbb{H}^n} \sum_{\ell = 1}^{k} c_{\ell}  (- \triangle_\Bbb{H})^{\ell} (\phi v_j) (\phi v_j) - \lambda (\phi v_j)^2 dV \\
  =& \int_{\Bbb{H}^n} \sum_{\ell = 1}^{k} c_{\ell} \sum_{ |\alpha| \geq 1}  D^\alpha \phi  D^{2\ell - \alpha} v_j (\phi v_j)  dV + \int_{\Bbb{H}^n} \sum_{\ell = 1}^{k} c_{\ell} \phi^2 v_j (- \triangle_\Bbb{H} )^{\ell} v_j dV - \lambda (\phi v_j)^2 dV \\
  =& \int_{\Bbb{H}^n} P_k v_j (\phi^2 v_j) - \lambda (\phi v_j)^2 dV + o(1).
\end{align*}

\noindent The limit is due to the fact that $\phi v_j$ is compactly supported and $v_j \rightharpoonup v = 0$ implies $L_{\text{loc}}^2$ strong convergence. Therefore,

$$
  < v_j , \phi^2 v_j >_\lambda = I_{\lambda,k} [\phi v_j] + o(1).
$$

\noindent By using H\"older's inequality, we obtain

\begin{align*}
  S_{\lambda,k} ||\phi v_j||_q^2 &\leq I_{\lambda,k} [\phi v_j] \\
    &= < v_j , \phi^2 v_j >_\lambda + o(1) \\
    &= \int_{\Bbb{H}^n} |v_j|^{q-2} (\phi v_j)^2 dV + o(1) \\
    &\leq ||\phi v_j||_q^2 \left( \int_{B((1,0), R)} |v_j|^q  dV\right)^{\frac{q-2}{q}} + o(1).
\end{align*}

\noindent This implies

$$
  \liminf_{n \to \infty}  \int_{B((1,0),R)} |v_j|^q  dV \geq S_{\lambda,k}^{\frac{q}{q-2}} \text{ and } \int_{\Bbb{H}^n \setminus B((1,0),R)} |v_j|^q dV \to 0.
$$

\noindent Therefore, from \eqref{P-S condition}, $I_{\lambda,k} [\phi v_j] \to S_{\lambda,k} $. Also, from our assumption, we know that $||\phi v_j||_2 \to 0$. Now
\begin{align*}
  S_{\lambda,k} &= \lim_{j \to \infty} I_{\lambda,k} [\phi v_j] \\
                    &= \lim_{j \to \infty} \frac{I_{\lambda,k} [\phi v_j]}{||\phi v_j||_q^2}  \\
                    &=  \lim_{j \to \infty} \frac{ \int_{\Bbb{H}^n} P_k (\phi v_j) (\phi v_j) dV }{||\phi v_j||_q^2} \\
                    &\geq  S_{n,k}.
\end{align*}
\noindent On the other hand, from Theorem \ref{Lu-Yang theorem}, $S_\lambda < S_{n,k}$, which is a contradiction. Hence the proof is completed.

\section{Nonexistence of $C^{2k}$ solutions: Proof of Theorem \ref{Nonexistence}}

In this section, we will establish the nonexistence of $C^{2k}$ solutions for the Brezis-Nirenberg problem on bounded star-shaped domains in the hyperbolic space $\mathbb{H}^n$, namely Theorem \ref{Nonexistence}.

\textbf{Proof of Theorem \ref{Nonexistence}}:  From \cite{Liu1, LuYang1},

  \begin{equation}\label{Transformation identity}
     \left( \frac{1 - |x|^2}{2} \right)^{k + \frac{n}{2}} (- \triangle)^k \left( \left( \frac{1 - |x|^2}{2} \right)^{k - \frac{n}{2}} u \right) = P_k u.
  \end{equation}
  Thus if we set $v = \left( \frac{1 - |x|^2}{2} \right)^{k - \frac{n}{2}} u$, we have the following equivalent equation on $\Omega' \subset \Bbb{R}^n$,

  $$
    (- \triangle)^k v = v^{q-1} + \lambda \left( \frac{2}{1 - |x|^2} \right)^{2k} v,
  $$
  where $\Omega'$ is also star-shaped. Moreover, since $\nabla_{\Bbb{H}^n}^\alpha v \vert_{\partial \Omega'} = 0$, we have $\nabla^\alpha v \vert_{\partial \Omega'} = 0$ for $|\alpha| \leq k - 1$ (this can be seen from the definition of $\nabla_{\Bbb{H}^n}$). For simplicity, we denote $p(x) =  \frac{2}{1 - |x|^2}$. We will try to establish a Pohozaev type identity by multiplying both sides of the equation by $x\nabla v$ and then use the  integration by parts. For the right hand side, we have

  \begin{align*}
    \int_{\Omega'} p^{2k} v x\cdot \nabla v dx &= \frac{1}{2} \int_{\Omega'} p^{2k} x\cdot \nabla v^2 dx, \\
                                                                      &= - \frac{1}{2} \int_{\Omega'} ( \nabla p^{2k} \cdot x + n p^{2k} ) v^2 dx,
  \end{align*}
  and

  \begin{align*}
    \int_{\Omega'} v^{q - 1} x\cdot  \nabla v dx
      &= \frac{1}{q} \int_{\Omega'} x \cdot \nabla v^q dx \\
      &= - \frac{n}{q} \int_{\Omega'} v^q dx.
  \end{align*}
  For the term $(- \triangle)^k v$, we first consider the case when $k$ is even,
  \begin{align*}
    &\int_{\Omega'} (- \triangle)^k v x \cdot \nabla v dx \\
      &= \int_{\Omega'} ( \triangle^{k/2 + 1} v ) \triangle^{k/2 - 1} (x \cdot \nabla v) dx \\
      &= \int_{\Omega'} ( \triangle^{k/2 + 1} v ) \triangle^{k/2 - 2} \nabla_j \nabla_j (x_i \nabla_i v) dx \\
      &= \int_{\Omega'} ( \triangle^{k/2 + 1} v ) \triangle^{k/2 - 2} \nabla_j ( \delta_{ji} \nabla_i v + x_i \nabla_j \nabla_i v ) dx \\
      &= \int_{\Omega'} ( \triangle^{k/2 + 1} v ) \triangle^{k/2 - 2} ( 2\triangle v + x_i \nabla_i \triangle v ) dx \\
      &= \int_{\Omega'} ( \triangle^{k/2 + 1} v ) \triangle^{k/2 - 3} (4\triangle^2 v + x_i \nabla_i \triangle^2 v ) dx \\
      &= \cdots \\
      &= \int_{\Omega'} ( \triangle^{k/2 + 1} v ) \left( (k -2) \triangle^{k/2 -1} v + x_i \nabla_i \triangle^{k/2 -1}v  \right) dx \\
      &= (k - 2) \int_{\Omega'} (\triangle^{k/2} v)^2 dx + \int_{\Omega'} \triangle^{k/2 + 1} v x \cdot \nabla (\triangle^{k/2 -1} v) dx \\
      &= (k - 2) \int_{\Omega'} (\triangle^{k/2} v)^2 dx  - \int_{\Omega'} \nabla_i \triangle^{k/2} v \nabla_i ( x_j \nabla_j \triangle^{k/2 -1} v ) dx \\
      &= (k - 2) \int_{\Omega'} (\triangle^{k/2} v)^2 dx - \int_{\Omega'} \nabla_i \triangle^{k/2} v ( \delta_{ij} \nabla_i \triangle^{k/2 -1} v + x_j \nabla_i \nabla_j \triangle^{k/2 -1} v ) dx \\
      &= (k - 2) \int_{\Omega'} (\triangle^{k/2} v)^2 dx + \int_{\Omega'} ( \triangle^{k/2} v )^2 dx - \int_{\partial \Omega'} \triangle^{k/2} v x_j \nu_i \nabla_j \nabla_j  \triangle^{k/2 -1} v \\
      &+ \int_{\Omega'} \triangle^{k/2} v ( \triangle^{k/2} v + x_j \nabla_j \triangle^{k/2} v ) dx\\
      &= k \int_{\Omega'} (\triangle^{k/2} v)^2 dx - \int_{\partial \Omega'} \triangle^{k/2} v \left( \frac{\partial^2}{\partial \nu^2} \triangle^{k/2 -1} v \right) (x \cdot \nu) d \sigma + \frac{1}{2} \int_{\Omega'} x_j \nabla_j (\triangle^{k/2} v)^2 dx \\
      &= k \int_{\Omega'} (\triangle^{k/2} v)^2 dx - \int_{\partial \Omega'} ( \triangle^{k/2} v )^2 (x \cdot \nu) d\sigma + \frac{1}{2} \int_{\partial \Omega'} ( \triangle^{k/2} v )^2 (x \cdot \nu) d\omega \\
      &- \frac{n}{2} \int_{\Omega'} (\triangle^{k/2} v)^2 dx \\
      &= - \frac{n - 2k}{2} \int_{\Omega'} (\triangle^{k/2} v)^2 dx - \frac{1}{2} \int_{\partial \Omega'} (\triangle^{k/2} v)^2 (x \cdot \nu) d\omega, \\
  \end{align*}
  where we make use of the assumption that $\nabla^\alpha v = 0$ on $\partial \Omega'$ when $|\alpha| \leq k -1$. When $k$ is odd, following a similar computation, we obtain

  $$
    \int_{\Omega'} (- \triangle)^k v x\nabla v dx = - \frac{n - 2k}{2} \int_{\Omega'} |\nabla \triangle^{\frac{k - 1}{2}} v|^2 dx - \frac{1}{2} \int_{\partial \Omega'} \left( \frac{\partial}{\partial \nu} \triangle^{\frac{k-1}{2}} v \right)^2 (x \cdot \nu) d\omega.
  $$
  Therefore when $k$ is even, we have

  \begin{align*}
     &\frac{n - 2k}{2} \int_{\Omega'} (\triangle^{k/2} v)^2 dx + \frac{1}{2} \int_{\partial \Omega'} (\triangle^{k/2} v)^2 (x \cdot \nu) d\omega \\
      &= \frac{\lambda}{2} \int_{\Omega'} ( \nabla p^{2k} \cdot x + n p^{2k} ) v^2 dx + \frac{n}{q} \int_{\Omega'} v^q dx.
  \end{align*}
  On the other hand, if we multiply $\frac{n-2k}{2} v$ to the original equation and integrate, we get

  $$
    \frac{n - 2k}{2} \int_{\Omega'} (\triangle^{k/2} v)^2 dx =  \lambda \frac{n - 2k}{2} \int_{\Omega'} p^{2k} v^2 dx +  \frac{n - 2k}{2} \int_{\Omega'} v^q dx.
  $$
  Therefore we have

  $$
    \frac{1}{2} \int_{\partial \Omega'} (\triangle^{k/2} v)^2 (x \cdot \nu) d\omega = \frac{\lambda}{2} \int_{\Omega'} ( \nabla p^{2k} \cdot x + n p^{2k} ) v^2 dx  - \lambda \frac{n - 2k}{2} \int_{\Omega'} p^{2k} v^2 dx.
  $$
Recalling the definition of $p(x)$, it is then easy to verify that

  $$
    \nabla p^{2k} \cdot x = 2k p^{2k - 1} \nabla p \cdot x = k p^{2k + 1} |x|^2.
  $$
  We hence obtain the Pohozaev type identity for the case when $k$ is even,

  $$
    \frac{1}{2} \int_{\partial \Omega'} (\triangle^{k/2} v)^2 (x \cdot \nu) d\omega = \frac{\lambda}{2} \int_{\Omega'} (p |x|^2 + 2) kp^{2k} v^2 dx.
  $$
  When $\lambda < 0$, apparently we have $v$ should be identically equal to zero and hence $u$ is trivial. For the case when $k$ is odd, following the same argument, we achieve the same conclusion. When $\lambda = 0$ and $k = 1$, notice that $x \cdot \nu > 0$ almost everywhere on $\partial \Omega'$, we have $\frac{\partial}{\partial \nu} v = 0$ almost everywhere. Now the equation

   $$
    (- \triangle)^k v = v^{q-1} + \lambda \left( \frac{2}{1 - |x|^2} \right)^{2k} v,
   $$
  gives us

  $$
    0 = \int_{\partial \Omega'} \frac{\partial}{\partial \nu} v d\sigma = - \int_{\Omega'} \triangle v dx = \int_{\Omega'} v^{q -1} dx.
  $$
  Therefore $v \equiv 0$ and hence $u \equiv 0$.
This completes the proof of Theorem \ref{Nonexistence}.

\section{Hyperbolic symmetry of the solution: Proofs of Theorem \ref{Integral equation theorem} and Theorem \ref{symmetry result}}

 We will first need  to give  the proof of Theorem \ref{Integral equation theorem}, namely
we  need to show that the positive solutions of the integral equation (\ref{Integral equation1}) are equivalent to those of the partial differential equation. Moreover, the integral kernel must satisfy certain conditions.
The proof of   Theorem \ref{Integral equation theorem} will be divided into two lemmas.

\begin{lemma}\label{integral equation lemma1}
   If $u \in W_0^{k,2}(\mathbb{H}^n)$ is a positive weak solution of the higher order differential equation
(\ref{differential equation}), then $u$ must satisfy the following integral equation for the Green function $G(x, y)$ of the operator $P_k - \lambda$
on the hyperbolic ball $\mathbb{B}^n$:
  \begin{equation}\label{Integral equation11}
     u(x) = \int_{\Bbb{B}^n} G(x,y) u^{q-1}(y) dV_y.
  \end{equation}

\end{lemma}

\begin{proof}
 Assume $u \in W_0^{k,2}(\mathbb{H}^n)$ is a positive weak solution of the higher order differential equation
  $$
    P_k u - \lambda u = u^{q - 1}.
  $$
 We start by recalling that the weak $W^{k,2}$ solution of the above differential equation is defined in following way: $u \in W^{1,2} (\Bbb{H}^n)$ is called the weak solution if for any $\phi \in C^\infty_0(\Bbb{H}^n)$, there holds

$$
  \int_{\Bbb{H}^n} \sum_{j = 1}^k c_j \nabla_\Bbb{H}^j u \nabla_\Bbb{H}^j \phi - \lambda u \phi dV = \int_{\Bbb{H}^n} u^{q-1} \phi dV.
$$
With the help of Fourier transform on $\Bbb{H}^n$ (See Section 2), such definition of weak solution is equivalent to:

$$
  D_n \int_{- \infty}^{+ \infty} \int_{\Bbb{S}^{n-1}} \left[ \prod_{j = 1}^k \left( \frac{\tau^2 + (2k - 1)^2}{4} \right)- \lambda \right] \hat{u}(\tau, \sigma) \hat{\phi}(\tau, \sigma) |\mathfrak{c}(\tau)|^{-2} d\sigma d\tau = \int_{\Bbb{H}^n} u^{q-1} \phi dV.
$$
Let $\psi $ satisfies $P_k \psi - \lambda \psi = \phi$, i.e. $\psi(x) = \int_{\Bbb{H}^n} G(x,y) \phi(y) dV_y$, which from a density argument, is still an applicable test function and under Fourier transform, it satisfies

$$
  \hat{\psi}(\tau, \sigma ) = \left[ \prod_{j = 1}^k \left( \frac{\tau^2 + (2k - 1)^2}{4} \right)- \lambda \right]^{-1} \hat{\phi}(\tau,\sigma).
$$
Now once we replace $\phi$ by $\psi$, we are to get

$$
  D_n \int_{- \infty}^{+ \infty} \int_{\Bbb{S}^{n-1}} \hat{u}(\tau, \sigma) \hat{\psi}(\tau, \sigma) |\mathfrak{c}(\tau)|^{-2} d\sigma d\tau = \int_{\Bbb{H}^n} u^{q-1}(x) \left( \int_{\Bbb{H}^n} G(x,y)\phi(y)dV_y \right) dV_x.
$$
Applying Plancherel formula on left hand side and changing order of integration on the right hand side, we get

$$
  \int_{\Bbb{H}^n} u(x) \phi(x) dV = \int_{\Bbb{H}^n} \left( \int_{\Bbb{H}^n} G(x,y) u^{q-1}(x) dV_x \right) \phi(y) dV_y
$$
holds for any $\phi \in C_0^\infty(\Bbb{H}^n)$, which instantly implies a solution of the differential equation is a solution of the integral equation. The other direction, i.e. the solution of the integral equation satisfies the differential equation, can be proved in the same spirit of applying Fourier transform and we omit it here.

\end{proof}

\begin{lemma}\label{integral equation lemma2} Let $G(x, y)$ be the Green function of the operator $P_k - \lambda$
on the hyperbolic ball $\mathbb{B}^n$.
Then, for any  fixed $y\in \mathbb{B}^n$,  $G(x,y)$ is a positive radially decreasing function with respect to the geodesic distance $\rho = d(x,y)$.
\end{lemma}

  \begin{proof} We start by recalling the definition of $P_k - \lambda$ for $0 < \lambda < \Lambda_1(P_k,\Bbb{H}^n) = \prod_{j = 1}^k \frac{(2j - 1)^2}{4}$, i.e.

\begin{align*}
  P_k - \lambda &= P_1 (P_1 + 2) \cdots (P_1 + k(k-1)) - \lambda \\
                &=(P_1 - \lambda_1) \cdots (P_1 - \lambda_l) (P_1^2 + A_1 P_1 + B_1) \cdots (P_1^2 + A_m P_1 + B_m ) \\
                &= (P_1 - \lambda_1) \cdots (P_1 - \lambda_l) \cdots (P_1 - \lambda_k),
\end{align*}
where $\lambda_1, \cdots,\lambda_k$ are roots (including complex roots) of the k-th order equation $x(x+2)\cdots(x+ k(k-1)) - \lambda = 0$. For those real $\lambda_j$'s, it is easy to see that $\lambda_j < \frac{1}{4}$, as otherwise we will have $\lambda_j (\lambda_j + 2) \cdots (\lambda_j + k(k-1)) \geq \prod_{j = 1}^k \frac{(2j - 1)^2}{4} > \lambda $. For those complex $\lambda_j$'s, we have $\textit{Re} \lambda_j < \frac{1}{4}$ since otherwise $\lambda = \lambda_j (\lambda_j + 2) \cdots (\lambda_j + k(k-1)) = |\lambda_j| |\lambda_j + 2| \cdots |\lambda_j + k(k-1)| \geq \textit{Re} \lambda_j (\textit{Re} \lambda_j + 2) \cdots (\textit{Re} \lambda_j + k(k-1))\geq \prod_{j = 1}^k \frac{(2j - 1)^2}{4} > \lambda$. Then we have the explicit formula of Green's function for each $(P_1 - \lambda_j)^{-1}$ and moreover, for each of the real valued $\lambda_j$'s, $(P_1 - \lambda_j)^{-1}$ is positive radially decreasing function (see \cite[Section 3]{LuYang1} for details).

We have to treat those complex valued $\lambda_j$'s carefully. Since they appear with their conjugates, without loss of generality, we only need to consider the function of the following form

$$
  (P_1 - \lambda_j)^{-1} * (P_1 - \overline{\lambda}_j)^{-1}.
$$
Since we have the formula of $(P_1 - \lambda_j)^{-1}$ in terms of the heat kernel:

$$
  (P_1 - \lambda_j)^{-1} = \int_0^\infty e^{t (\triangle_{\Bbb{H}^n} + \frac{n(n-2)}{4} + \lambda_j)} dt,
$$
we have

\begin{align*}
  &(P_1 - \lambda_j)^{-1} * (P_1 - \overline{\lambda}_j)^{-1} \\
  =& \left( \int_0^\infty e^{t \triangle_{\Bbb{H}^n}} e^{t (\frac{n(n-2)}{4} + \lambda_j)} dt \right) * \left( \int_0^\infty e^{s \triangle_{\Bbb{H}^n}} e^{s (\frac{n(n-2)}{4} + \overline{\lambda}_j)} ds \right) \\
  =& \int_{[0, \infty) \times [0,\infty)} e^{\frac{n(n-2)}{4} (t + s)} e^{t \lambda_j + s \overline{\lambda}_j} e^{t \triangle_{\Bbb{H}^n}} * e^{s\triangle_{\Bbb{H}^n} } ds dt \\
  =& \int_{[0, \infty) \times [0,\infty)} e^{\frac{n(n-2)}{4} (t + s)} e^{s \lambda_j + t \overline{\lambda}_j} e^{t \triangle_{\Bbb{H}^n}} * e^{s\triangle_{\Bbb{H}^n} } ds dt \\
  =& \overline{(P_1 - \lambda_j)^{-1} * (P_1 - \overline{\lambda}_j)^{-1}}.
\end{align*}
This means $(P_1 - \lambda_j)^{-1} * (P_1 - \overline{\lambda}_j)^{-1}$ is a real function and from the monotonicity of the heat kernel, it is also easy to see that $(P_1 - \lambda_j)^{-1} * (P_1 - \overline{\lambda}_j)^{-1}$ is positive. Recall that when $n = 2m + 1$,  $e^{t \triangle_{\Bbb{H}^n}}$ is given explicitly by the formula (see Davies \cite{Davies}):

$$
  e^{t \triangle_{\Bbb{H}^n}} = 2^{-m - 2} \pi^{-m - 1/2} t^{-1/2} e^{- \frac{(n-1)^2}{4}t} \left( - \frac{1}{\sinh \rho} \frac{\partial}{\partial \rho} \right)^m e^{-\frac{\rho^2}{4t}}.
$$
Combining this formula with the semigroup property of $e^{t \triangle_{\Bbb{H}^n}}$, we have

\begin{align*}
  & \frac{\partial}{\partial \rho} (P_1 - \lambda_j)^{-1} * (P_1 - \overline{\lambda}_j)^{-1} \\
   =& (-\sinh \rho) \left( -\frac{1}{\sinh \rho} \frac{\partial}{\partial \rho}\right) (P_1 - \lambda_j)^{-1} * (P_1 - \overline{\lambda}_j)^{-1} \\
   =& - 2^{-m - 2} \pi^{-m - 1/2} \sinh \rho \int_{[0, \infty) \times [0,\infty)} (t + s)^{-1/2} e^{- \frac{t + s}{4}} e^{t \lambda_j + s \overline{\lambda}_j} \left( - \frac{1}{\sinh \rho} \frac{\partial}{\partial \rho} \right)^{m+1} e^{-\frac{\rho^2}{4t}}  ds dt  \\
   =& - C \sinh \rho (\tilde{P}_1 - \lambda_j)^{-1} * (\tilde{P}_1 - \overline{\lambda}_j)^{-1},
\end{align*}
where $(\tilde{P}_1 - \lambda_j)^{-1}$ stands for the kernel function for dimension $\tilde{n} = 2m + 3$ and $C$ is some positive constant. This instantly implies that $(P_1 - \lambda_j)^{-1} * (P_1 - \overline{\lambda}_j)^{-1}$ is decreasing.

On the other hand, when $n = 2m$, $e^{t \triangle_{\Bbb{H}^n}}$ is given explicitly by the formula (see \cite{Davies}):

$$
  e^{t \triangle_{\Bbb{H}^n}} = (2 \pi)^{- \frac{n+1}{2}} t^{-1/2} e^{- \frac{(n-1)^2}{4} t} \int_\rho^\infty \frac{\sinh r}{\sqrt{\cosh r - \cosh \rho }} \left( - \frac{1}{\sinh r} \frac{\partial}{\partial r} \right)^m e^{-\frac{r^2}{4t}} dr.
$$
Then we have

\begin{align*}
  &(P_1 - \lambda_j)^{-1} * (P_1 - \overline{\lambda}_j)^{-1} \\
  =& (2 \pi)^{- \frac{n+1}{2}} \int_{[0, \infty) \times [0,\infty)} (t + s)^{-1/2} e^{-\frac{t + s}{4} } e^{t \lambda_j + s \overline{\lambda}_j} \left( \int_\rho^\infty \frac{\sinh r}{\sqrt{\cosh r - \cosh \rho }} \left( - \frac{1}{\sinh r} \frac{\partial}{\partial r} \right)^m e^{-\frac{r^2}{4(t+s)}} dr \right) ds dt \\
  =& - 2 (2 \pi)^{- \frac{n+1}{2}} \int_{[0, \infty) \times [0,\infty)} (t + s)^{-1/2} e^{-\frac{t + s}{4} } e^{t \lambda_j + s \overline{\lambda}_j} \cdot\\ & \left( \int_\rho^\infty \sqrt{\cosh r - \cosh \rho } \frac{\partial }{\partial r} \left( - \frac{1}{\sinh r} \frac{\partial}{\partial r} \right)^m e^{-\frac{r^2}{4(t+s)}} dr \right) ds dt.
\end{align*}
Hence

\begin{align*}
  & \frac{\partial}{\partial \rho} (P_1 - \lambda_j)^{-1} * (P_1 - \overline{\lambda}_j)^{-1} \\
  =&  (2 \pi)^{- \frac{n+1}{2}}  \int_{[0, \infty) \times [0,\infty)} (t + s)^{-1/2} e^{-\frac{t + s}{4} } e^{t \lambda_j + s \overline{\lambda}_j}  \\
  & \left(  \int_\rho^\infty \frac{\sinh \rho}{\sqrt{\cosh r - \cosh \rho }} \frac{\partial}{\partial r} \left( - \frac{1}{\sinh r} \frac{\partial}{\partial r} \right)^m e^{-\frac{r^2}{4(t+s)}} dr \right) ds dt \\
  =& -(2 \pi)^{- \frac{n+1}{2}}  \int_{[0, \infty) \times [0,\infty)} (t + s)^{-1/2} e^{-\frac{t + s}{4} } e^{t \lambda_j + s \overline{\lambda}_j} \int_\rho^\infty \frac{\sinh \rho \sinh r}{\sqrt{\cosh r - \cosh \rho} }
  \left( - \frac{1}{\sinh r} \frac{\partial}{\partial r} \right)^{m + 1} \\
   &e^{-\frac{r^2}{4(t+s)}} dr  ds dt \\
  =& - C \sinh \rho (\tilde{P}_1 - \lambda_j)^{-1} * (\tilde{P}_1 - \overline{\lambda}_j)^{-1},
\end{align*}
where $(\tilde{P}_1 - \lambda_j)^{-1}$ stands for the kernel function for dimension $\tilde{n} = 2m + 2$ and $C$ is some positive constant. Therefore $(P_1 - \lambda_j)^{-1} * (P_1 - \overline{\lambda}_j)^{-1} $ is positive and decreasing.

Recall that $G(x, y)$ is the Green function of the operator $P_k-\lambda$ on the hyperbolic ball $\mathbb{B}^n$.
    Given any fixed $y\in \mathbb{B}^n$, we claim that $G(x,y)$ is a positive radially decreasing function with respect to the geodesic distance $\rho = d(x,y)$. Since the Green function $G(x, y)$ of the operator $P_k - \lambda = (P_1 - \lambda_1) \cdots (P_1 - \lambda_l) \cdots (P_1 - \lambda_k)$ is the convolution of positive decreasing functions, to show that $G(x, y)$ is positive decreasing, we first {\bf claim } that for any positive functions $H_1(x,y), H_2(x,y) : (\Bbb{B}^n \times \Bbb{B}^n) \setminus \{ x = y \} \to \Bbb{R}$, if both of them are decreasing with respect to $d(x,y)$, then $L(x,y) = \int_{\Bbb{H}^n} H_1(x,z)H_2(z,y) dV_z$ is also decreasing with respect to $d(x,y)$.

     To prove this claim, we first show $L(x,y) = L(d(x,y))$. This is because for any isometry $T: \Bbb{B}^n \to \Bbb{B}^n$, we have
  \begin{align*}
    L(Tx,Ty) &= \int_{\Bbb{B}^n} H_1(Tx,z)H_2(z,Ty) dV_z \\
                  &= \int_{\Bbb{B}^n} H_1(x,T^{-1} z) H_2(T^{-1}z, y) dV_z \\
                  &= \int_{\Bbb{B}^n} H_1(x,z)H_2(z,y) dV_z \\
                  &= L(x,y).
  \end{align*}
  Now fixing $y$, consider any geodesic ray from $y$ and any two points $x,\overline{x}$ in this ray. Without loss of generality, we can assume $d(x,y) \leq d(\overline{x},y)$ and it suffices to show $L(x,y) \geq L(\overline{x}, y)$. For simplicity, denote $I = \{ z \in \Bbb{B}^n : d(z,x) \leq d(z,\overline{x}) \}$ and $II = \Bbb{B}^n \setminus I$. From the foliation structure of $\Bbb{B}^n$, there is a unique reflection $I_t$ on $\Bbb{B}^n$ such that $I_t(x) = I_t (\overline{x})$.
  We further denote $\overline{z} = I_t(z)$. then
  \begin{align*}
    &L(x,y) - L(\overline{x},y) \\
    &=  \int_{\Bbb{B}^n} H_1(x,z)H_2(z,y) dV_z -   \int_{\Bbb{B}^n} H_1(\overline{x},z)H_2(z,y) dV_z \\
    &= \left( \int_I H_1(x,z)H_2(z,y) dV_z + \int_{II} H_1(x,z)H_2(z,y) dV_z \right) \\
    &-   \left( \int_I H_1(\overline{x},z)H_2(z,y) dV_z + \int_{II} H_1(\overline{x},z)H_2(z,y) dV_z \right)\\
    &= \left( \int_I H_1(x,z)H_2(z,y) dV_z  + \int_{I} H_1(x,\overline{z})H_2(\overline{z},y) dzV_z \right) \\
    &-  \left( \int_I H_1(\overline{x},z)H_2(z,y) dV_z  + \int_{I} H_1(\overline{x},\overline{z})H_2(\overline{z},y) dzV_z \right) \\
    &= \int_{I} (H_1(x,z) - H_1(\overline{x},z)) H_2(z,y) dV_z - \int_{I} (H_1(\overline{x},\overline{z}) - H_1(x,\overline{z})) H_2(\overline{z},y) dV_z \\
    &= \int_I (H_1(x,z) - H_1(\overline{x},z)) H_2(z,y) dV_z - \int_I (H_1(x,z) - H_1(\overline{x},z)) H_2(\overline{z},y) dV_z \\
    &= \int_I (H_1(x,z) - H_1(\overline{x},z)) (H_2(z,y) - H_2(\overline{z} , y)) dV_y.
  \end{align*}
  Since $d(x,z) \leq d(\overline{x},z)$ and $d(z,y) \leq d(\overline{z},y)$, we have $L(x,y) \geq L(\overline{x},y)$ and our claim is therefore proved. This completes the proof of Lemma \ref{integral equation lemma2}.
\end{proof}
Combining Lemma \ref{integral equation lemma1} and Lemma \ref{integral equation lemma2},  this finishes the proof of
   Theorem \ref{Integral equation theorem}. \\

We are now ready to prove the symmetry result, i.e., Theorem  \ref{symmetry result}.

{\bf Proof of Theorem  \ref{symmetry result}.} Because of the validity of  Theorem \ref{Integral equation theorem}, it suffices to prove
the symmetry of the solution to the integral equation (\ref{Integral equation1}).

  To perform the moving plane argument, we fix one specific direction and for such direction, consider $t > 0$, $U_t$ splits the hyperbolic ball $\Bbb{B}^n$ into two parts. We denote $\Sigma_t = \cup_{s < t} U_s$. For any $x \in \Sigma_t$, denote $\overline{t} = I_t (x)$ and $u_t (x) = u(\overline{t})$. We have

    \begin{align*}
    &u(x) - u_t(x) \\
    =& \int_{\Sigma_t} G(x,y)u(y)^{q-1} dV_y + \int_{\Sigma_t^c} G(x,y)u(y)^{q-1} dV_y \\
    -& \int_{\Sigma_t} G(x,y)u_t(y)^{q-1} dV_y - \int_{\Sigma_t^c} G(x,y)u_t(y)^{q-1} dV_y\\
    =& \int_{\Sigma_t} G(x,y)u(y)^{q-1} dV_y + \int_{\Sigma_t} G(x,\overline{y})u(\overline{y})^{q-1} dV_y \\
    -&\int_{\Sigma_t} G(x,y)u_t(y)^{q-1} dV_y -\int_{\Sigma_t} G(x,\overline{y})u_t(\overline{y})^{q-1} dV_y \\
    =& \int_{\Sigma_t} G(x,y)u(y)^{q-1} dV_y + \int_{\Sigma_t} G(\overline{x}, y)u_t(y)^{q-1} dV_y \\
    -&\int_{\Sigma_t} G(x,y)u_t(y)^{q-1} dV_y -\int_{\Sigma_t} G(\overline{x},y)u(y)^{q-1} dV_y \\
    =& \int_{\Sigma_t} \left( G(x,y) - G(\overline{x},y) \right) \left( u(u)^{q-1} - u_t(y)^{q-1} \right) dV_y,
  \end{align*}
    where we use the fact that $G(x,\overline{y}) = G(\overline{x},y)$. This is because $I_t$ is isometry and hence $d(x,\overline{y}) = d(\overline{x},y)$.

  \medskip

  One of the key steps is to compare $G(x,y)$ and $G(\overline{x},y)$. This is equivalent to compare $d(x,y)$ and $d(\overline{x},y)$. We claim that $d(x,y) < d(\overline{x},y)$. It suffices to prove that $d(A_{-t}(x), A_{-t}(y)) < d(A_{-t}(\overline{x}), A_{-t}(y))$. Let $\gamma_1$, $\gamma_2$ and $\gamma_3$ be geodesic curves connecting $A_{-t}(x)$ and $A_{-t}(\overline{x})$, $A_{-t}(x)$ and $A_{-t}(y)$, $A_{-t}(\overline{x})$ and $A_{-t}(y)$. Also, let $\gamma_4$ be the geodesic curve passing through $A_{-t}(y)$ and perpendicular to $\gamma_1$ and we denote the joint point as $Q$. Now we have two right hyperbolic triangles: $\triangle A_{-t}(x) Q A_{-t}(y)$ and $\triangle A_{-t}(\overline{x}) Q A_{-t}(y)$. Since $d(A_{-t}(x) , Q) < d( A_{-t}(\overline{x}), Q)$, from basic hyperbolic geometry, we conclude that $d(A_{-t}(x), A_{-t}(y) ) < d (A_{-t}(\overline{x}), A_{-t}(y))$ and hence $d(x,y) < d(\overline{x},y)$. Therefore, we have $G(x,y) > G(\overline{x},y)$.

  \medskip

  In $\Sigma_t$, we denote $\Sigma_t^- = \{ x \in \Sigma_t: u_t (x) > u(x) \}$, then for $x \in \Sigma_t^-$, we have from \eqref{Green function estimate},

  \begin{align*}
    &u_t (x) - u (x) \\
    &\leq C \int_{\Sigma_t^-} G(x,y) u_t^{q-2}(y) \left( u_t(y) - u(y) \right) dV_y \\
    &\leq C \int_{\Sigma_t^-}  \left( \frac{1}{\sinh \frac{d(x,y)}{2}} \right)^{n-2k} u_t^{q-2}(y) \left( u_t(y) - u(y) \right) dV_y
  \end{align*}
  From the Hardy-Littlewood-Sobolev inequality on $\Bbb{H}^n$ (Theorem \ref{Hardy-Littlewood-Sobolev on H^n}) and H\"older's inequality, we have

  \begin{align*}
    || u_t - u ||_{L^q (\Sigma_t^-)}
      &\leq C \left( \int_{\Sigma_t^-} u_t(y)^q dV_y \right)^{2k/n} || u_t - u ||_{L^q (\Sigma_t^-)} \\
      &\leq C \left( \int_{\Sigma_t^c} u(y)^q dV_y \right)^{2k/n} || u_t - u ||_{L^q (\Sigma_t^-)}
  \end{align*}
  Since $u \in L^q (\Bbb{H}^n)$, we can hence pick $t$ large enough such that $C \left( \int_{\Sigma_t^c} u(y)^q dV_y \right)^{2k/n} < 1$. This implies for $t$ large enough, $|| u_t - u ||_{L^q (\Sigma_t^-)}  = 0$ and hence $\Sigma_t^-$ is of measure zero.

  \medskip

  Now we shift $U_t$ as long as $u \geq u_t$ holds in $\Sigma_t$. Suppose for $t_0$, we have $u(x) \geq u_{t_0} (x)$ on $\Sigma_{t_0}$ but $u \neq u_t$. For $t_0 - \epsilon < t \leq t_0$ we can again prove

  $$
    || u_t - u ||_{L^q (\Sigma_{t}^-)}
      \leq C \left( \int_{I_t (\Sigma_{t}^-)} u(y)^q dV_y \right)^{2k/n} || u_t - u ||_{L^q (\Sigma_t^-)} .
  $$
  When $\epsilon$ is small enough, $\Sigma_t^-$ is close enough to zero such that $C \left( \int_{I_t (\Sigma_{t}^-)} u(y)^q dV_y \right)^{2k/n} < 1$. Then $|| u_t - u ||_{L^q (\Sigma_t^-)} = 0$ which means we can keep on shifting $U_{t_0}$.
This completes the proof of  Theorem \ref{symmetry result}.
  \medskip

\section{A revisit of Theorem \ref{Existence result at critical exponent on bounded domain} in lower dimensions: Proof of Theorem \ref{Lower bound lambda}}

In the second statement of Theorem \ref{Existence result at critical exponent on bounded domain}, we give a sufficient condition to ensure the existence of solution for Problem \ref{Brezis-Nirenberg on bounded domain} when $2k + 1 \leq n \leq 4k-1$. In this section, we will particularly consider the case $k = 2$ (so that $n = 5,6,7$) and $\Omega = B_R(0)$ is the hyperbolic ball centered at the origin with radius $R$ for some $0 < R < 1$. We aim to find a wider range for $\lambda$ so that Problem \ref{Brezis-Nirenberg on bounded domain} still has a solution. Such refinement relies on finding the ``optimal'' choice of a smooth cut-off function and one of the key ingredients is the hyperbolic symmetry of the solution, which helps to reduce the problem to ODE's.

We first recall some well-known results of Legendre functions. Associated Legendre functions $P_l^{m}(\cosh \rho)$ and $P_l^{-m}(\cosh \rho)$ are solutions to the following Legendre equation:

\begin{equation}\label{Legendre equation}
  u''(\rho) + \coth \rho u'(\rho) + \left( l(l+1) - \frac{m^2}{\sinh^2 \rho} \right) u(\rho) = 0.
\end{equation}
The associated Legendre function has the following expression:

$$
  P_l^m (\cosh \rho) = \frac{1}{\Gamma (1 - m)} \coth^m (\rho/2) F_{2,1}\left[ l,l+1,1-m;-\sinh^2 (\rho/2) \right],
$$
where for complex numbers $a,b,c,z$, the hypergeometric function is given by

$$
  F_{2,1} [a,b,c;z] = \sum_{n = 0}^\infty \frac{(a)_n (b)_n}{(c)_n} \frac{z^n}{n!},
$$
where $(\beta)_n = \prod_{j = 0}^{n-1} (\beta + j)$, for any $\beta \in \Bbb{C}$. Moreover, the associated Legendre functions satisfy the following raising and lowering relations:

$$
  \frac{d}{d \rho} P_l^m(\cosh \rho) = \frac{1}{\sinh \rho} P_l^{m+1}(\cosh \rho) + \frac{m \cosh \rho}{\sinh^2 \rho} P_l^m (\cosh \rho)
$$
and

$$
 \frac{d}{d \rho} P_l^{m+1}(\cosh \rho) = \frac{l(l+1) - m(m+1)}{\sinh \rho} P_l^m (\cosh \rho) - \frac{(m+1) \cosh \rho}{\sinh^2 \rho} P_l^{m+1} (\cosh \rho).
$$

Based on Theorem \ref{symmetry result}, we are searching for non-trivial radially symmetric solutions of

$$
  \begin{cases}
    P_2 u(\rho) - \lambda u(\rho) = |u(\rho)|^{\frac{8}{n-4}}u(\rho) \ &\textit{on } B_R(0), \\
    u = \frac{\partial u}{\partial \nu} = 0 \ &\textit{on } \partial B_R(0),
  \end{cases}
$$
where $\rho$ denotes the geodesic distance on $\Bbb{H}^n$. Recall that when $k = 2$, the identity \eqref{Transformation identity} becomes

$$
   \left( \frac{1 - |x|^2}{2} \right)^{2 + \frac{n}{2}} (- \triangle)^2 \left( \left( \frac{1 - |x|^2}{2} \right)^{2 - \frac{n}{2}} u \right) = P_2 u
$$
and if we set $v = \left( \frac{1 - |x|^2}{2} \right)^{2 - \frac{n}{2}} u$,

\begin{align*}
  S_{\lambda,2}
    &\leq \frac{\int_\Omega P_2 u \cdot u - \lambda u^2 dV}{\left( \int_\Omega |u|^{\frac{2n}{n-4}} dV \right)^{\frac{n-4}{n}}} \\
    &= \frac{\int_\Omega \left[ (-\triangle)^2 v \cdot v \left( \frac{1 - |x|^2}{2} \right)^n - \lambda v^2 \left( \frac{1 - |x|^2}{2} \right)^{n-4} \right] \left( \frac{2}{1 - |x|^2} \right)^n dx }{  \left(\int_\Omega |v|^{\frac{2n}{n-4}} \left( \frac{1 - |x|^2}{2} \right)^{\frac{n-4}{2} \cdot \frac{2n}{n-4}} \left( \frac{2}{1 - |x|^2} \right)^n dx \right)^{\frac{n-4}{n}}} \\
    &= \frac{\int_\Omega (-\triangle)^2 v \cdot v - \lambda \left( \frac{2}{1 - |x|^2} \right)^4 v^2 dx}{\left( \int_\Omega |v|^{\frac{2n}{n-4}} dx \right)^{\frac{n-4}{n}}}.
\end{align*}
Now we let $\phi(|x|)= \phi(r)$ be a nonnegative smooth cut-off function defined on $\Omega=B_R(0)$ satisfying $\phi = \phi' = 0$ on $\partial B_R(0)$. If we replace $v$ in the above quotient by

$$
  v_\epsilon = \frac{\phi(r)}{(\epsilon^2 + r^2)^{\frac{n-4}{2}}},
$$
then following a quite similar calculation as \cite[Lemma 4,1]{EdmundsFortunatoJannelli1}, we obtain a sufficient condition for $S_{\lambda,2} < S_{n,2}$ to hold, which further implies the existence of non-trivial solution. This condition reads as follows:

\begin{prop}
If there exists a non-negative function $\phi(r) \in C^\infty (B_R(0))$ satisfying $\phi = \phi' = 0$ on $\partial B_R(0)$ and

$$
  \int_0^R |\phi''|^2 r^{7-n} dr + 3(n-3) \int_0^R |\phi'|^2 r^{5-n} dr - \lambda \int_0^R |\phi|^2 \left( \frac{2}{1-r^2} \right)^4 r^{7-n} dr < 0,
$$
then $S_{\lambda,2} < S_{n,2}$.
\end{prop}
From the criterion we established in the proof of Theorem \ref{Existence result at critical exponent on bounded domain}, there  exists a non-trivial solution as long as

$$
  \lambda > \lambda^* = \inf_{\phi \in C^\infty(B_R(0)), \phi \geq 0, \phi(R) = \phi'(R) = 0} \frac{\int_0^R |\phi''|^2 r^{7-n} dr + 3(n-3) \int_0^R |\phi'|^2 r^{5-n} dr}{\int_0^R |\phi|^2 \left( \frac{2}{1-r^2} \right)^4 r^{7-n} dr}.
$$
The corresponding Euler equation for $\phi$ is

$$
  \phi^{(4)} + 2(7-n)\frac{\phi'''}{r} + (n^2 - 16n + 51)\frac{\phi''}{r^2} + 3(n-5)(n-3) \frac{\phi'}{r^3} - \lambda^* \left( \frac{2}{1-r^2} \right)^4 \phi = 0.
$$
Let $\Phi(r) = r^{4-n} \phi(r)$, the above ODE can be rewritten as

$$
  \Phi^{(4)} + \frac{2(n-1)}{r}\Phi''' + (n^2 - 4n + 3)\left( \frac{\Phi''}{r^2} - \frac{\Phi'}{r^3}  \right) - \lambda^* \left( \frac{2}{1-r^2} \right)^4 \Phi = 0,
$$
which is exactly

$$
  (-\triangle)^2 \Phi - \lambda^* \left( \frac{2}{1 - r^2} \right)^6 \Phi = 0.
$$
Now let $\Phi (r) = \left( \frac{1 - r^2}{2} \right)^{2 - \frac{n}{2}} \Psi (r)$ and apply \eqref{Transformation identity} again, we get

\begin{align*}
  0 &= P_2 \Psi - \lambda^* \Psi \\
    &=(P_1^2 + 2P_1 - \lambda^*) \Psi \\
    &= (P_1 - \lambda_1^*)(P_1 - \lambda_2^*)\Psi,
\end{align*}
where $\lambda_1^* = -1 + \sqrt{1 + \lambda^*}$ and $\lambda_2^* = -1 - \sqrt{1 + \lambda^*}$ are solutions of the quadratic equation $y^2 + 2y - \lambda^* = 0$. From classical ODE theory, $\Psi$ is a linear combination of the solution of

$$
  (P_1 - \lambda_j^*)\Psi = 0, \ j = 1,2.
$$
Recall the definition of $P_1 = - \triangle_{\Bbb{H}^n} - \frac{n(n-2)}{4}$ and $\triangle_{\Bbb{H}^n} = -\frac{\partial^2}{\partial^2 \rho} - (n - 1) \coth \rho \frac{\partial}{\partial \rho}$ when it acts on radial functions, the above equation can hence be written as

$$
  \frac{\partial^2}{\partial^2 \rho} \Psi + (n - 1) \coth \rho \frac{\partial}{\partial \rho} \Psi + \left( \lambda_j^* + \frac{n (n-2)}{4} \right) \Psi = 0.
$$
Now we further set $\Psi (\rho) = \sinh^{\frac{2-n}{2}} \rho \psi(\rho)$, then the above equation is equivalent to the following standard Legendre equation

$$
  \psi'' + \coth \rho \psi' + \left(\lambda_j^* - \frac{(\frac{n-2}{2})^2}{\sinh^2 \rho}\right) \psi = 0.
$$
Therefore $\psi (\rho)$ is a linear combination of $P_{l_1^*}^{\pm \frac{n-2}{2}} (\cosh \rho)$. Now combining the information of $\phi, \Phi, \Psi, \psi$, a straightforward simplification gives

$$
  \phi = \cosh^{\frac{-n+2}{2}} \frac{\rho}{2} \sinh^{\frac{n-6}{2}} \frac{\rho}{2} \left( C_1 P_{\lambda_1^*}^{\frac{2-n}{2}} + C_2 P_{\lambda_1^*}^{\frac{n-2}{2}} + C_3 P_{\lambda_2^*}^{\frac{2-n}{2}} + C_4 P_{\lambda_2^*}^{\frac{n-2}{2}} \right)
$$
for some  $C_1,\, C_2,\, C_3,\, C_4 \in \Bbb{R}$.
Since we only look for smooth $\phi$'s, from the expression of the Legendre function, $C_2 = C_4 = 0$. In addition, since we require $\phi = \phi' = 0$ on $\partial B_R(0)$, if we denote $\overline{\rho} = \log \frac{1 + R}{1 - R}$, we have

$$
  C_1 P_{\lambda_1^*}^{\frac{2-n}{2}} (\overline{\rho}) + C_3 P_{\lambda_2^*}^{\frac{2-n}{2}} (\overline{\rho}) = 0
$$
and

$$
  C_1 (P_{\lambda_1^*}^{\frac{2-n}{2}})'(\overline{\rho}) + C_3 (P_{\lambda_2^*}^{\frac{2-n}{2}})' (\overline{\rho}) = 0.
$$
In order to achieve a non-trivial solution, we need

$$
\det
  \begin{pmatrix}
    P_{\lambda_1^*}^{\frac{2-n}{2}} (\overline{\rho}) & P_{\lambda_2^*}^{\frac{2-n}{2}} (\overline{\rho}) \\
    (P_{\lambda_1^*}^{\frac{2-n}{2}})' (\overline{\rho}) & (P_{\lambda_2^*}^{\frac{2-n}{2}})' (\overline{\rho})
  \end{pmatrix}
  = 0.
$$
By using the raising and lowering relations of the Legendre function, after some basic linear algebra simplification, we have

$$
\det
  \begin{pmatrix}
    P_{\lambda_1^*}^{\frac{2-n}{2}} (\overline{\rho}) & P_{\lambda_2^*}^{\frac{2-n}{2}} (\overline{\rho}) \\
    P_{\lambda_1^*}^{\frac{4-n}{2}} (\overline{\rho}) & P_{\lambda_2^*}^{\frac{4-n}{2}} (\overline{\rho})
  \end{pmatrix}
  = 0.
$$
Therefore we have the following refined theorem:

\begin{theorem}
  When $k = 2$ and consider $4 < n < 8$ and $\Omega = B_R(0)$. If we define $\lambda^*$ as the first positive value such that

  $$
    \det
  \begin{pmatrix}
    P_{\lambda_1^*}^{\frac{2-n}{2}} (\overline{\rho}) & P_{\lambda_2^*}^{\frac{2-n}{2}} (\overline{\rho}) \\
    P_{\lambda_1^*}^{\frac{4-n}{2}} (\overline{\rho}) & P_{\lambda_2^*}^{\frac{4-n}{2}} (\overline{\rho})
  \end{pmatrix}
  = 0,
  $$
  Then Problem \ref{Brezis-Nirenberg on bounded domain} has a  non-trivial solution when $\lambda^* < \lambda < \overline{\Lambda} (P_2 , \Omega)$
\end{theorem}

\begin{remark}
 One implicit fact adapted in Theorem \ref{Lower bound lambda} is that the constant $\lambda^*$ satisfies $\lambda^* \leq \overline{\Lambda}(P_2, \Omega)$. To see this, recall in the definition of $\lambda^*$, we have

 $$
   \lambda^* \leq \frac{\int_0^R |\phi''|^2 r^{7-n} dr + 3(n-3) \int_0^R |\phi'|^2 r^{5-n} dr}{\int_0^R |\phi|^2 \left( \frac{2}{1-r^2} \right)^4 r^{7-n} dr}
 $$
 for any $\phi \in C^\infty(B_R(0))$ satisfying $\phi(R) = \phi'(R) = 0$. Now if we define $\Phi = r^{4-n} \phi$ and $\tilde{\Phi} = \left( \frac{1 - |x|^2}{2} \right)^{-2 + \frac{n}{2}} \Phi$, we have the following

 \begin{align*}
   &\frac{\int_0^R |\phi''|^2 r^{7-n} dr + 3(n-3) \int_0^R |\phi'|^2 r^{5-n} dr}{\int_0^R |\phi|^2 \left( \frac{2}{1-r^2} \right)^4 r^{7-n} dr} \\
   =& \frac{\int_0^R |\Phi''|^2 r^{n-1} dr + (n-1) \int_0^R |\Phi'|^2 r^{n-3} dr}{\int_0^R |\Phi|^2 \left( \frac{2}{1-r^2} \right)^4 r^{n-1} dr} \\
   =& \frac{\int_{B_R(0)} (-\triangle)^2 \Phi \cdot \Phi dx}{\int_{B_R(0)} |\Phi|^2 \left( \frac{2}{1 - |x|^2} \right)^4 dx} \\
   =& \frac{\int_\Omega P_2 \tilde{\Phi} \cdot \tilde{\Phi} dV}{\int_\Omega |\tilde{\Phi}|^2 dV},
 \end{align*}
 where we have used the relation \eqref{Transformation identity} for the last line. In particular, when $\tilde{\Phi}$ is chosen to be the eigenfunction of $P_2$ corresponding to the first eigenvalue $\overline{\Lambda}(P_2,\Omega)$, it directly give us the validity of $\lambda^* \leq \overline{\Lambda}(P_2,\Omega)$.
\end{remark}

{\bf Acknowledgement.} The main results of this paper have been presented by the first author at ``Workshop on harmonic analysis and applications" at Tsinghua Sanya International Mathematics Forum in January, 2020, at the online AMS special session ``Geometric Inequalities and Nonlinear Partial Differential Equations" in September, 2020 and the online seminar ''Geometric and functional inequalities and applications" at reserachseminars.org in September, 2020 and the online AMS special session on ``Geometric and Functional Inequalities and Nonlinear Partial Differential Equations" in March, 2021.

\bibliographystyle{plain}

\end{document}